\documentclass[10.5pt]{article}
\usepackage{amsmath}
\usepackage{amssymb}
\usepackage{amscd}
\usepackage{xspace}
\usepackage{verbatim}
\newcommand{\mysection}[1]{
\section{#1}\setcounter{equation}{0}}
\title{\bf Boundary singularities of solutions of $N$-harmonic 
equations with absorption
\footnote {To appear in{ \it Journal of Functional Analysis}}
}%% 
\author{{\bf Rouba Borghol}, {\bf Laurent V\'eron}\\
{\small Department of Mathematics,}\\
 {\small  University of Tours,  FRANCE}
}%%

\date{}
\begin{document}
\maketitle
{\small\it Dedicated to Jim Serrin, on his eightieth birthday}\\

\noindent{\small {\bf Abstract} We study the boundary behaviour of 
solutions $u$ of $-\Delta_{N}u+ |u|^{q-1}u=0$ in a bounded smooth 
domain $\Omega\subset\mathbb R^{N}$ subject to the boundary condition $u=0$ except at one 
point, in the range $q>N-1$. We prove that if $q\geq 2N-1$ such a $u$ 
is identically zero, while, if $N-1<q<2N-1$, $u$ inherits a boundary 
behaviour which either corresponds to a weak singularity, or to  a strong 
singularity. Such singularities are effectively constructed.} 
\vspace{1mm}
\hspace{.05in}

%% FONT commands
\newcommand{\txt}[1]{\;\text{ #1 }\;}%% Used in math only
\newcommand{\tbf}{\textbf}%% Bold face. Usage: \tbf{...}
\newcommand{\tit}{\textit}%% Italic
\newcommand{\tsc}{\textsc}%% Small caps
\newcommand{\trm}{\textrm}
\newcommand{\mbf}{\mathbf}%% Math bold
\newcommand{\mrm}{\mathrm}%% Math Roman
\newcommand{\bsym}{\boldsymbol}%% Bold math symbol
%%Macros for changing font size in math.
\newcommand{\scs}{\scriptstyle}%% as in subscript
\newcommand{\sss}{\scriptscriptstyle}%% as in sub-subscript
\newcommand{\txts}{\textstyle}
\newcommand{\dsps}{\displaystyle}
%%Macros for changing font size in text.
\newcommand{\fnz}{\footnotesize}
\newcommand{\scz}{\scriptsize}
%%\tiny<\scz<\fsz<\small<\large<\Large<\huge<\Huge
%%%%%%%%%%%%
%%%%%%%%%%%%
%% EQUATION commands
\newcommand{\be}{
\begin{equation}
}
\newcommand{\bel}[1]{
\begin{equation}
\label{#1}}
\newcommand{\ee}{
\end{equation}
}%% This macro does not work with amstex.
\newcommand{\eqnl}[2]{
\begin{equation}
\label{#1}{#2}
\end{equation}
}%%use not advisable; confusing
%%%%%%%%%%%%%%%
%% Unnumbered THEOREM env.
%% New env. to be used for unnumbered theorem, lemma etc. (but with specified name)
\newtheorem{subn}{\name}
\renewcommand{\thesubn}{}
\newcommand{\bsn}[1]{\def\name{#1}
\begin{subn}}
\newcommand{\esn}{
\end{subn}}
%%%%%%%%%%%%%%
%% NUMBERED THEOREM env.
%% Environments: theorem, lemma, corollary defintion and related commands,
%% designed to provide consecutive numbering of these forms.
\newtheorem{sub}{\name}[section]
\newcommand{\dn}[1]{\def\name{#1}}   %used in conjuction with sub or subn.
\newcommand{\bs}{
\begin{sub}}
\newcommand{\es}{
\end{sub}}
\newcommand{\bsl}[1]{
\begin{sub}\label{#1}}
%% the above must be preceeded by \dn (name definition),
%% however this is superceded by the list of commands bth etc.  below.
%%%%%%%%%%%%
%% NUMBERED THEOREM env. (cont.)
%% List of commands derived from 'sub' env. for theorem, lemma etc.
%% designed to provide consecutive numbering of these forms.
\newcommand{\bth}[1]{\def\name{Theorem}
\begin{sub}\label{t:#1}}
\newcommand{\blemma}[1]{\def\name{Lemma}
\begin{sub}\label{l:#1}}
\newcommand{\bcor}[1]{\def\name{Corollary}
\begin{sub}\label{c:#1}}
\newcommand{\bdef}[1]{\def\name{Definition}
\begin{sub}\label{d:#1}}
\newcommand{\bprop}[1]{\def\name{Proposition}
\begin{sub}\label{p:#1}}
%%%%%%%%%%%%%%%%%%%%%%%%%%%%%%%%%%
%% RERERENCE commands.
%% \newcommand{\R}[1]{(\ref{#1})}
\newcommand{\R}{\eqref}
\newcommand{\rth}[1]{Theorem~\ref{t:#1}}
\newcommand{\rlemma}[1]{Lemma~\ref{l:#1}}
\newcommand{\rcor}[1]{Corollary~\ref{c:#1}}
\newcommand{\rdef}[1]{Definition~\ref{d:#1}}
\newcommand{\rprop}[1]{Proposition~\ref{p:#1}}
%%%%%%%%%%%
%% ARRAY commands.
\newcommand{\BA}{
\begin{array}}
\newcommand{\EA}{
\end{array}}
\newcommand{\BAN}{\renewcommand{\arraystretch}{1.2}
\setlength{\arraycolsep}{2pt}
\begin{array}}
\newcommand{\BAV}[2]{\renewcommand{\arraystretch}{#1}
\setlength{\arraycolsep}{#2}
\begin{array}}
%Note: The first variable gives the amount of stretching: (#1) x default.
%For instance #1=1.2 means a 20% stretching. The second variable should be
%written for instance in the form  4pt ; here the default is 5pt
%\newcommand{\EAN}{\end{array}\setlength{\arraycolsep}{5pt}}
\newcommand{\BSA}{
\begin{subarray}}
\newcommand{\ESA}{
\end{subarray}}
%Note: These are used in subscripts as well as superscripts. They work essentially
%% like 'array'.
\newcommand{\BAL}{
\begin{aligned}}
\newcommand{\EAL}{
\end{aligned}}
\newcommand{\BALG}{
\begin{alignat}}
\newcommand{\EALG}{
\end{alignat}}%% the abbrev. does not work with latex2e
\newcommand{\BALGN}{
\begin{alignat*}}
\newcommand{\EALGN}{
\end{alignat*}}%% the abbrev. does not work with latex2e
%% The 'aligned' environment must be placed inside an 'equation' env.
%% in the same way as the array.
%% One could use also the 'align' env. or the 'alignat' env.
%% However in this case each line is numbered, unless '\notag' is used.
%% The 'alignat'
%% has a slightly different format (the number of columns must be specified in advance)
%% but it has the advantage that the distance between columns is at our disposition.
%% (The default would be zero distance.) Using 'alignat*' we can have the advantages
%% of alignat plus the situation where separate lines are not numbered.
%% However in this case there is no numbering at all (unless we provide a tag).
%%%%%%%%%%
%% PROOF, REMARK etc.
\newcommand{\note}[1]{\textit{#1.}\hspace{2mm}}
\newcommand{\Proof}{\note{Proof}}
\newcommand{\qeda}{\hspace{10mm}\hfill $\square$}
\newcommand{\qed}{\\
${}$ \hfill $\square$}
\newcommand{\Remark}{\note{Remark}}
%%%%%%%% Style command.
\newcommand{\modin}{$\,$\\
[-4mm] \indent}
%% To be used after \mysection in order to start new line with \indent.
%%%%%%%%%%%%
%% MATHEMATICAL symbols
\newcommand{\forevery}{\quad \forall}
\newcommand{\set}[1]{\{#1\}}
\newcommand{\setdef}[2]{\{\,#1:\,#2\,\}}
\newcommand{\setm}[2]{\{\,#1\mid #2\,\}}
%% Arrows
\newcommand{\lra}{\longrightarrow}
\newcommand{\lla}{\longleftarrow}
\newcommand{\llra}{\longleftrightarrow}
\newcommand{\Lra}{\Longrightarrow}
\newcommand{\Lla}{\Longleftarrow}
\newcommand{\Llra}{\Longleftrightarrow}
\newcommand{\warrow}{\rightharpoonup}
%% Brackets, delimiters
\newcommand{
\paran}[1]{\left (#1 \right )}%% adjustable parantheses
\newcommand{\sqbr}[1]{\left [#1 \right ]}%% adjustable square brackets
\newcommand{\curlybr}[1]{\left \{#1 \right \}}%% adjustable curly brackets
\newcommand{\abs}[1]{\left |#1\right |}%% adjustable vertical delimiters
\newcommand{\norm}[1]{\left \|#1\right \|}%% adjustable norm
\newcommand{
\paranb}[1]{\big (#1 \big )}%% non-adjustable parantheses (big)
\newcommand{\lsqbrb}[1]{\big [#1 \big ]}%% non-adjustable square brackets (big)
\newcommand{\lcurlybrb}[1]{\big \{#1 \big \}}%% non-adjustable curly brackets (big)
\newcommand{\absb}[1]{\big |#1\big |}%% non-adjustable vertical delimiters (big)
\newcommand{\normb}[1]{\big \|#1\big \|}%% non-adjustable norm (big)
\newcommand{
\paranB}[1]{\Big (#1 \Big )}%% non-adjustable parantheses (Big)
\newcommand{\absB}[1]{\Big |#1\Big |}%% non-adjustable vertical delimiters (Big)
\newcommand{\normB}[1]{\Big \|#1\Big \|}%% non-adjustable norm (Big)

%%%%%%%%%%%%%%%%%
%% Adjustable parantheses etc. in a different DEFINITION format.
%\def\adp(#1){\left (#1 \right )}%% adjustable parantheses
%\def\adsb(#1){\left [#1\right ]}%% adjustable square brackets
%\def\adcb(#1){\left \{#1\right \}}%% adjustable curly brackets
%\def\abs|#1|{\left |#1\right |}%% adjustable vertical delimiters
%%%%%%%%%%%%%%%%
%% More mathematical symbols
\newcommand{\thkl}{\rule[-.5mm]{.3mm}{3mm}}
\newcommand{\thknorm}[1]{\thkl #1 \thkl\,}
\newcommand{\trinorm}[1]{|\!|\!| #1 |\!|\!|\,}
\newcommand{\bang}[1]{\langle #1 \rangle}%% angle bracket
\def\angb<#1>{\langle #1 \rangle}%% angle bracket
%% The two last lines yield the same result.
%% The second is used as follows: \angb<a,b>
\newcommand{\vstrut}[1]{\rule{0mm}{#1}}
\newcommand{\rec}[1]{\frac{1}{#1}}
%% OPERATOR names.
%% OPERATOR names.
\newcommand{\opname}[1]{\mbox{\rm #1}\,}
\newcommand{\supp}{\opname{supp}}
\newcommand{\dist}{\opname{dist}}
\newcommand{\myfrac}[2]{{\displaystyle \frac{#1}{#2} }}
\newcommand{\myint}[2]{{\displaystyle \int_{#1}^{#2}}}
\newcommand{\mysum}[2]{{\displaystyle \sum_{#1}^{#2}}}
\newcommand {\dint}{{\displaystyle \int\!\!\int}}%%%%%%%%%%
%%%%%%% SPACE commands
\newcommand{\q}{\quad}
\newcommand{\qq}{\qquad}
\newcommand{\hsp}[1]{\hspace{#1mm}}
\newcommand{\vsp}[1]{\vspace{#1mm}}
%%%%%%%%%%%
%% ABREVIATIONS
\newcommand{\ity}{\infty}
\newcommand{\prt}{
\partial}
\newcommand{\sms}{\setminus}
\newcommand{\ems}{\emptyset}
\newcommand{\ti}{\times}
\newcommand{\pr}{^\prime}
\newcommand{\ppr}{^{\prime\prime}}
\newcommand{\tl}{\tilde}
\newcommand{\sbs}{\subset}
\newcommand{\sbeq}{\subseteq}
\newcommand{\nind}{\noindent}
\newcommand{\ind}{\indent}
\newcommand{\ovl}{\overline}
\newcommand{\unl}{\underline}
\newcommand{\nin}{\not\in}
\newcommand{\pfrac}[2]{\genfrac{(}{)}{}{}{#1}{#2}}% frac with parantheses.
%%%%%%%%%%%
%%%%%%%%%%%%%

%%Macros for Greek letters.
\def\ga{\alpha}     \def\gb{\beta}       \def\gg{\gamma}
\def\gc{\chi}       \def\gd{\delta}      \def\ge{\epsilon}
\def\gth{\theta}                         \def\vge{\varepsilon}
\def\gf{\phi}       \def\vgf{\varphi}    \def\gh{\eta}
\def\gi{\iota}      \def\gk{\kappa}      \def\gl{\lambda}
\def\gm{\mu}        \def\gn{\nu}         \def\gp{\pi}
\def\vgp{\varpi}    \def\gr{\rho}        \def\vgr{\varrho}
\def\gs{\sigma}     \def\vgs{\varsigma}  \def\gt{\tau}
\def\gu{\upsilon}   \def\gv{\vartheta}   \def\gw{\omega}
\def\gx{\xi}        \def\gy{\psi}        \def\gz{\zeta}
\def\Gg{\Gamma}     \def\Gd{\Delta}      \def\Gf{\Phi}
\def\Gth{\Theta}
\def\Gl{\Lambda}    \def\Gs{\Sigma}      \def\Gp{\Pi}
\def\Gw{\Omega}     \def\Gx{\Xi}         \def\Gy{\Psi}

%%Macros for calligraphic letters.
\def\CS{{\mathcal S}}   \def\CM{{\mathcal M}}   \def\CN{{\mathcal N}}
\def\CR{{\mathcal R}}   \def\CO{{\mathcal O}}   \def\CP{{\mathcal P}}
\def\CA{{\mathcal A}}   \def\CB{{\mathcal B}}   \def\CC{{\mathcal C}}
\def\CD{{\mathcal D}}   \def\CE{{\mathcal E}}   \def\CF{{\mathcal F}}
\def\CG{{\mathcal G}}   \def\CH{{\mathcal H}}   \def\CI{{\mathcal I}}
\def\CJ{{\mathcal J}}   \def\CK{{\mathcal K}}   \def\CL{{\mathcal L}}
\def\CT{{\mathcal T}}   \def\CU{{\mathcal U}}   \def\CV{{\mathcal V}}
\def\CZ{{\mathcal Z}}   \def\CX{{\mathcal X}}   \def\CY{{\mathcal Y}}
\def\CW{{\mathcal W}} \def\CQ{{\mathcal Q}} 
%%%%%
%%Macros for 'blackboard' letters (See (27) for display.)
\def\BBA {\mathbb A}   \def\BBb {\mathbb B}    \def\BBC {\mathbb C}
\def\BBD {\mathbb D}   \def\BBE {\mathbb E}    \def\BBF {\mathbb F}
\def\BBG {\mathbb G}   \def\BBH {\mathbb H}    \def\BBI {\mathbb I}
\def\BBJ {\mathbb J}   \def\BBK {\mathbb K}    \def\BBL {\mathbb L}
\def\BBM {\mathbb M}   \def\BBN {\mathbb N}    \def\BBO {\mathbb O}
\def\BBP {\mathbb P}   \def\BBR {\mathbb R}    \def\BBS {\mathbb S}
\def\BBT {\mathbb T}   \def\BBU {\mathbb U}    \def\BBV {\mathbb V}
\def\BBW {\mathbb W}   \def\BBX {\mathbb X}    \def\BBY {\mathbb Y}
\def\BBZ {\mathbb Z}

%%Macros for Ghotic (Fraktur) letters.
\def\GTA {\mathfrak A}   \def\GTB {\mathfrak B}    \def\GTC {\mathfrak C}
\def\GTD {\mathfrak D}   \def\GTE {\mathfrak E}    \def\GTF {\mathfrak F}
\def\GTG {\mathfrak G}   \def\GTH {\mathfrak H}    \def\GTI {\mathfrak I}
\def\GTJ {\mathfrak J}   \def\GTK {\mathfrak K}    \def\GTL {\mathfrak L}
\def\GTM {\mathfrak M}   \def\GTN {\mathfrak N}    \def\GTO {\mathfrak O}
\def\GTP {\mathfrak P}   \def\GTR {\mathfrak R}    \def\GTS {\mathfrak S}
\def\GTT {\mathfrak T}   \def\GTU {\mathfrak U}    \def\GTV {\mathfrak V}
\def\GTW {\mathfrak W}   \def\GTX {\mathfrak X}    \def\GTY {\mathfrak Y}
\def\GTZ {\mathfrak Z}   \def\GTQ {\mathfrak Q}

\font\Sym= msam10 % special symbols
\def\SYM#1{\hbox{\Sym #1}}
\newcommand{\bdw}{\prt\Gw\xspace}
\medskip
\mysection {Introduction}
Let $\Gw$ be a domain is $\BBR^N$ ($N\geq 2$) with a $C^2$ compact boundary 
$\prt\Gw$. Let $g$ be a continous real valued function and $a\in \prt\Gw$. 
This paper deals 
with the study of solutions 
$u\in C^{1}(\bar\Gw\setminus\{a\})$ of the problem
\begin {equation}\label {main1}\left\{\BA {l}
-\,div\left(\abs {Du}^{N-2}Du\right)+g(u)=0\qquad\mbox {in } \Gw\\
\phantom {-\div\left(\abs {Du}^{N-2}Du\right)_{N}u+,, }
u=0\qquad\mbox {on } \prt\Gw\setminus \{a\},
\EA\right.\end {equation}
and we shall be more specifically interested in the case when $g$ has a 
power growth at infinity. When $N=2$, this problem falls into the 
scope of the boundary singularity problem for semilinear elliptic 
equations. The study of the $N$-dimensional problem
\begin {equation}\label {lin1}\left\{\BA {l}
-\Gd u+g(u)=0\qquad\mbox {in } \Gw\\
\phantom {-----}
u=0\qquad\mbox {on } \prt\Gw\setminus \{a\},
\EA\right.\end {equation}
has been initiated by Gmira and V\'eron in \cite {GV}. Among the 
subjects under consideration were the question of removability of 
isolated boundary singularities and, in the case such singularities 
do exist, their precise description. This seminal article was at the origin 
of a long series of further works by Dynkin, Kuznetsov, Le Gall, 
Marcus and V\'eron in the framework of the trace theory and, later 
on, the fine trace theory in the case where $g(r)=r\abs r^{q-1}$, $q>1$. 
One of the main reasons for such a large 
impact consists of the observation of the existence of a critical 
exponent $q=q^{*}=(N+1)/(N-1)$. If $q\geq q^{*}$ any solution of
\begin {equation}\label {lin2}\left\{\BA {l}
-\Gd u+\abs u^{q-1}u=0\qquad\mbox {in } \Gw\\
\phantom {-----,,-}\!
u=0\qquad\mbox {on } \prt\Gw\setminus \{a\},
\EA\right.\end {equation}
is identically zero, while if $1<q<q^{*}$ it appears that there 
exist two possible behaviours of singular solutions near $a$, the 
solutions with weak singularities and the ones with the strong singular 
behaviour. Later on, these two types of singular solutions played a 
fundamental role in the description of the rough trace of positive 
solutions of (\ref {lin2}). \smallskip

Although the techniques needed are considerably more refined, it 
appeared that the description of solutions of (\ref {main1}) inherits the 
same structure as for (\ref {lin1}). The first step is to understand 
the model case problem
\begin {equation}\label {main2}\left\{\BA {l}
-\,div\left(\abs {Du}^{N-2}Du\right)+\abs u^{q-1}u=0\qquad\mbox {in } \Gw\\
\phantom {-\div\left(\abs {Du}^{N-2}Du\right)_{N}u+,;;,,,}
u=0\qquad\mbox {on } \prt\Gw\setminus \{a\},
\EA\right.\end {equation}
To this equation, we associate the homogeneous equation
\begin {equation}\label {hom2}\left\{\BA {l}
-\,div\left(\abs {Du}^{N-2}Du\right)=0\qquad\mbox {in } \Gw\\
\phantom {-\div\left(\abs {Du}^{N-2}Du\right)}
u=0\qquad\mbox {on } \prt\Gw\setminus \{a\}.
\EA\right.\end {equation}
It is proved in \cite {BV} that for any $k>0$ there exists a unique 
solution $u=u_{k}$ of (\ref{hom2}) satisfying 
\begin {equation}\label {weak1}
u_{k}(x)=k\myfrac {\gr(x)}{\abs {x-a}^2}(1+\circ (1))\quad \mbox {as 
}x\to a,\;(x-a)/\abs {x-a}\to \gs,
\end {equation}
where $\gr(x)=\dist (x,\prt\Gw)$. When $k=1$, this solution plays the role of the Poisson kernel, although neither any weak formulation nor any reasonable trace theory seems to exists, and we shall denote it by $V_a^\Gw$. The behaviour (\ref{weak1}) (up to a multiplicative constant) corresponds to {\it weak singularity behaviour} 
for (\ref {main1}), whenever such singularities exist.  The first 
result we prove is the following:\smallskip

\nind{\bf Theorem }{\it Let $N-1<q<2N-1:=q_{c}$ . Then for any $k>0$ there exists a unique 
solution $u=u_{k,a}$ of problem (\ref {main2}) satisfying
\begin {equation}\label {weak2}
u_{k,a}(x)=k\myfrac {\gr(x)}{\abs {x-a}^2}(1+\circ (1))\quad \mbox {as 
}x\to a,\;(x-a)/\abs {x-a}\to \gs.
\end {equation}
Furthermore $u_{\infty,a}=\lim_{k\to\infty}$ 
exists and is a solution of (\ref {main2}) which satisfies
\begin {equation}\label {strong1}
\lim_{\scriptsize \BA {c}x\to a\\ \frac{x-a}{\abs {x-a}}\to\gs\EA}
\abs {x-a}^{N/(q+1-N)}u_{\infty,a}(x)=\gw(\gs),
\end {equation}
and $\gw$ 
is the unique positive solution of the following quasilinear 
equation on the upper hemisphere  of the unit sphere $S^{N-1}$,
\begin{equation} \label{strong1'}\left\{\BA {l}
-div_{\gs}\left(\left(\gb^2_{q}\gw^2+\abs {\nabla_{\gs}\gw}^{2}\right)^{(N-2)/2}
\nabla_{\gs}\gw\right)\\
\phantom{---------;}
-\Gl\left(\gb^2_{q}\gw^2+\abs {\nabla_{\gs}\gw}^{2}
\right)^{(N-2)/2}\!\!\!\!\!\!\gw+\abs \gw^{q-1}\gw
=0
\;\;\mbox { on }S^{N-1}_{+}\\[2mm]
\phantom{--------------------------;}
\gw=0\;\;\mbox {on }\;\prt S^{N-1}_{+},
\EA\right.
\end{equation}
}\\
where $\gb_q=N/(q+1-N)$ and $\Gl=(N-1)\gb_q^2$. 
The proof of the existence of $u_{k,a}$, as well as its singular behaviour, is settled upon the conformal invariance of the $N$-harmonic operator and the construction of subsolution of the same equation. Estimate (\ref{strong1}) is proved by scaling method. The role of the critical exponent 
$q_c=2N-1$ is enlighted by the following result.\medskip 

%%%%%%%%%%%%%%%%%%%%%%%%%%%%%%%%%%%%%%%%%%%%%%%%%%%%%%%%%%%%%%%%%%%%%%%%%%%%%%%%%%%%%%%%
%%%%%%THEOREM 1%%%%%%%%%%%%%%%%%%%%%%%%%%%%%%%%%%%%%%%%%%%%%%%%%%%%%%%%%%%%%%%%%%%%%%%%%%%
\nind{\bf Theorem }{\it Let $g$ be a 
continuous function such that 
\begin {equation}\label {remov1}\BA {l}
(i)\quad\liminf_{r\to\infty} g(r)/r^{q_{c}}>0\\[2mm]
(ii)\quad\limsup_{r\to-\infty} g(r)/\abs{r}^{q_{c}}<0.
\EA\end {equation}
Then any function $u\in C^{1}(\overline\Gw\setminus \{a\})$ solution of 
(\ref {main1}) extends as a function $\tilde u\in C(\overline\Gw)$.}\\

As in the semilinear case, the occurrence coincides with the case where the blow-up exponent $-\gb_q$ which is natural for equation (\ref{main2}) coincides with the one of the function $V^\Gw_a$ solution of (\ref{hom2}). Finally we provide the full classification of 
positive solutions of 
problem (\ref {main2}).\smallskip 

\nind{\bf Theorem }{\it Let $N-1<q<q_{c}$ and $u$ is any nonnegative solution of (\ref 
{main2}), then\smallskip

\nind (i) Either $u\equiv 0$,\smallskip

\nind (ii) Either there  exists $k>0$ such that $u=u_{k,a}$. \smallskip

\nind (iii) Or $u=u_{\infty,a}$.}\\

In the proof of (iii) the boundary Harnack inequalities that 
satisfies any positive solution of (\ref{main2}) (see \cite{BBV}) play a fundamental role.\\

Our paper is organized as follows\smallskip

\nind 1- Introduction\smallskip

\nind 2- Weak and strong boundary singularities\smallskip

\nind 3- The removability result\smallskip

\nind 4- The classification theorem

%%%%%%%%%%%%%%%%%%%%%%%%%%%%%%%%%%%%%%%%%%%%%%%%%%%%%%%%%%%%%%%%%%%%%%%%%%%%%%%%%%%%%%%%
%%%%%%SECTION 2 %%%%%%%%%%%%%%%%%%%%%%%%%%%%%%%%%%%%%%%%%%%%%%%%%%%%%%%%%%%%%%%%%%%%%%%%%%%
\mysection {Weak and strong boundary singularities}
The construction of positive solutions of 
\begin {equation}\label {equ+abs}
-div\left(\abs {Du}^{N-2}Du\right)+\abs u^{q-1}u=0,
\end{equation}
is settled upon three facts: the existence of solutions to the 
homogeneous equation
\begin {equation}\label {equ+hom}
-div\left(\abs {Du}^{N-2}Du\right)=0,
\end{equation}
the conformal invariance of (\ref {equ+hom}) and an a priori estimate 
 satisfied by any solution of \!(\ref{equ+abs}). Throughout this paper 
 $C$ denotes a positive constant which depends only on the structural 
 assumptions corresponding to $N$, $p$, $q$ and $\Gw$. The value of 
 the constant may change from one occurrence to another. 
 %%%%%%%%%%%%%%%%%%%%%%%%%%%%%%%%%%%%%%%%%%%%%%%%%%%%%%%%%%%%%%%%%%%%%%%%%%%%%%%%%%%%%%%%
%%%%PROPOSITION   KELLER-OSSERMAN I %%%%%%%%%%%%%%%%%%%%%%%%%%%%%%%%%%%%%%%%%%%%%%%%%%%%%%%%%%%%%%
%%%%%%%%%%%%%%%%%%%%%%%%%%%%%%%%%%%%%%%%%%%%%%%
\bprop {KOV} Let $\Gw\subset\BBR^{N}$ be a domain with a compact 
boundary and $a\in\prt\Gw$. Consider real numbers $q> p-1>0$, $A>0$ 
and $B\geq 0$. If $u\in C(\overline\Gw\setminus \{a\})\cap 
W^{1,p}_{loc}(\Gw)$ is a weak solution of 
\begin {equation}\label {sub}\left\{\BA {l}
-div\left(\abs {Du}^{p-2}Du\right)+A\abs u^{q-1}u\leq B\quad \mbox 
{in }\Gw\\\phantom{--------------}
u\leq 0\quad \mbox {on }\prt\Gw\setminus\{a\},
\EA\right.\end{equation}
it satisfies
\begin {equation}\label {up-est}
u(x)\leq \left(\myfrac {\gl}{A\abs{x-a}^p}\right)^{1/(q+1-p)}+
\left(\myfrac {\gm B}{A}\right)^{1/q}\quad\forevery 
x\in\overline\Gw\setminus\{a\},
\end{equation}
where $\gl$ and $\gm$ depends on $N$, $p$ and $q$.
\es
\Proof By assumption
\begin {equation}\label{weak}
\myint{\Gw}{}\left(\abs {Du}^{p-2}Du.D\gz+A\abs u^{q-1}u\gz\right)dx
\leq B\myint{\Gw}{}\gz dx
\end {equation}
for any $\gz\in W^{1,p}(\Gw)$ with compact support, $\gz\geq 0$. Let 
$\eta\in C^{2}(\BBR)$ be a nonnegative function such that 
$0\leq\eta'\leq 1$, $\eta''\geq 0$, $\eta=\eta'=\eta''$ on 
$(-\infty,0]$, $0<\eta(r)\leq r$ on $(0,\infty)$. For $\ge>0$ we set 
$\eta_{\ge}(r)=\eta ((r-\ge)_{+})$. Let $\gz\in 
W^{1,p}(\BBR^{N}\setminus \{0\})$ 
with compact support. Inasmuch $(\eta'_{\ge}(u))^{p-1}\gz$ has compact 
support in $\Gw$ and
$$D\left((\eta'_{\ge}(u))^{p-1}\gz\right)=
(\eta'_{\ge}(u))^{p-1}D\gz 
+(p-1)(\eta'_{\ge}(u))^{p-2}\eta''_{\ge}(u)\gz Du,$$
it belongs to $W^{1,p}(\Gw)$ and is an admissible test function 
for (\ref{weak}). Thus
$$\BA {l}
\myint{\Gw}{}\left(\abs 
{Du}^{p-2}Du.D\left((\eta'_{\ge}(u))^{p-1}\gz\right)
+A\abs u^{q-1}u(\eta'_{\ge}(u))^{p-1}\gz\right)dx
\leq B\myint{\Gw}{}(\eta'_{\ge}(u))^{p-1}\gz dx,
\EA$$
and
$$\abs 
{Du}^{p-2}Du.D\left((\eta'_{\ge}(u))^{p-1}\gz\right)\geq 
(\eta'_{\ge}(u))^{p-1}\abs{Du}^{p-2}Du.D\gz 
=\abs{Dv_{\ge}}^{p-2}Dv_{\ge}.D\gz,
$$
where we have set $v_{\ge}=\eta_{\ge}(u)$. Furthermore, $\eta$ can be 
chosen such that $r^q(\eta_{\ge}'(r))^{p-1}\geq \eta_{\ge}^q(r)$, for 
example if we fix $\eta(r)=r^{2}/2\gd$ on $(0,\gd]$ and 
$\eta(r)=r-\gd/2$ on $[\gd,\infty)$ for some $\gd>0$. We extend 
$v_{\ge}$ by $0$ outside $\overline\Gw\setminus\{a\}$ and denote by 
$\tilde v_{\ge}$ the 
new function, then $\tilde v_{\ge}\in 
W^{1,p}_{loc}(\BBR^{N}\setminus\{a\})\cap C(\BBR^{N}\setminus\{a\})$ 
and
\begin {equation}\label{weak2'}
\myint{\Gw}{}\left(\abs {D\tilde v_{\ge}}^{p-2}D\tilde v_{\ge}.D\gz
+A\abs {\tilde v_{\ge}}^{q-1}\tilde v_{\ge}\gz\right)dx
\leq B\myint{\Gw}{}\gz dx.
\end {equation}
This means that $\tilde v_{\ge}$ is a weak subsolution in 
$\BBR^{N}\setminus\{a\}$. By \cite [Lemma 1.3]{VV}, we derive 
$$\tilde v_{\ge}(x)\leq \left(\myfrac {\gl}{A\abs{x-a}^p}\right)^{1/q+1-p}+
\left(\myfrac {\gm B}{A}\right)^{1/q}\quad\forevery 
x\in\BBR^{N}\setminus\{a\},
$$
for some $\gl>0$ and $\gm>0$ depending on $N$, $p$ and $q$. Letting 
successively $\ge\to 0$ and $\gd\to 0$ we obtain (\ref{sub}).
\qeda \\

When $\Gw$ is smooth we have a sharper estimate
 %%%%%%%%%%%%%%%%%%%%%%%%%%%%%%%%%%%%%%%%%%%%%%%%%%%%%%%%%%%%%%%%%%%%%%%%%%%%%%%%%%%%%%%%
%%%%PROPOSITION KELLER-OSSERMAN II %%%%%%%%%%%%%%%%%%%%%%%%%%%%%%%%%%%%%%%%%%%%%%%%%%%%%%%%%%%%%%%
%%%%%%%%%%%%%%%%%%%%%%%%%%%%%%%%%%%%%%%%%%%%
\bprop {UPREG} Let $\Gw\subset\BBR^{N}$ be a bounded domain with $C^{2}$
boundary and $a\in\prt\Gw$. Let $q\geq p-1>1$ and $a>0$. 
If $u\in C(\overline\Gw\setminus \{a\})\cap 
W^{1,p}_{loc}(\Gw)$ is a weak solution of (\ref{sub}) with $B=0$, 
there exists $C>0$ depending on $\Gw$, $p$ and $q$ such that
\begin {equation}\label {up-est2}
u(x)\leq \myfrac {C\gr(x)}{\left(A\abs{x-a}^{q+1}\right)^{1/(q+1-p)}}
\forevery x\in\overline\Gw\setminus\{a\},
\end{equation}
where $\gr(x)=\dist (x,\prt\Gw$.
\es
\Proof By translation we can assume that $a=0$. For $\ge>0$ let 
$v_{\ge}$ be the solution of
\begin {equation}\label {subepsilon}\left\{\BA {l}
-div\left(\abs {Dv_{\ge}}^{p-2}Dv_{\ge}\right)+A\abs v_{\ge}^{q-1}v_{\ge}=0,
\quad \mbox 
{in }\Gw^{\ge}=\Gw\setminus B_{\ge}\\\phantom{--------------;}
v_{\ge}= u_{+}\quad \mbox {on }\prt\Gw^{\ge}.
\EA\right.\end{equation}
By \cite [Lemma 1.3]{VV} as in the proof of \rprop {KOV} and the maximum 
principle, there holds
$$u_{+}(x)\leq v_{\ge}(x)\leq \left(\myfrac {\gl}{A(\abs{x}-\ge)^p}\right)^{1/(q+1-p)}
\quad\forevery 
x\in\overline\Gw^\ge.
$$
Consequently $\ge\leq \ge'\Longrightarrow v_{\ge}\geq v_{\ge'}$. 
Letting $\ge\to 0$ and using the previous inequalities and the 
classical regularity results for solutions of quasilinear equations 
\cite {Lie} we conclude that $v_{\ge}$ converges, as $\ge\to 0$, to 
some $v$ which is a nonnegative solution of 
\begin {equation}\label {sub0}\left\{\BA {l}
-div\left(\abs {Dv}^{p-2}Dv\right)+A v^{q}=0,\quad \mbox 
{in }\Gw\\\phantom{-----------}
v= 0\quad \mbox {on }\prt\Gw\setminus \{0\}.
\EA\right.\end{equation}
and dominate $u$. Further, if $\ell>0$  
the function $v^{\ell}$ defined by 
$v^{\ell}(y)=\ell^{p/(q+1-p)}v(\ell y)$ is a solution of (\ref{sub0}) 
with $\Gw$ replaced by $\Gw_{\ell}=\ell^{-1}\Gw$. Let 
$x\in\overline\Gw\setminus \{0\}$ and $\ell=\abs x$. Since 
$$0\leq v^{\ell}(y)\leq \left(\myfrac 
{\gl}{A(\abs{y})^p}\right)^{1/(q+1-p)}\forevery y\in \Gw_{\ell},
$$
and
$$\max \left\{\abs{Dv^\ell(y)}:y\in \Gw_{\ell}\cap 
B_{3/2}\setminus B_{2/3}\right\}
\leq M\max \left\{\abs {v^\ell(z)}:z\in \Gw_{\ell}\cap 
B_{2}\setminus B_{1/2}\right\},
$$
where $M$ is uniformly bounded because the curvature of 
$\prt\Gw_{\ell}$ is bounded, we obtain that $Dv^\ell(y)$ is uniformly 
bounded by some constant $C$ on $\Gw_{\ell}\cap B_{3/2}\setminus B_{2/3}$. 
Because 
$Dv^\ell(y)=\ell^{(q+1)/(q+1-p)}Dv(\ell y)$, it follows that
$$\abs {Dv(x)}\leq \myfrac {C }{A^{1/q+1-p}\abs x^{(q+1)/(q+1-p)}}
$$
By the mean value Theorem, and using the fact that $v$ vanishes on 
$\prt\Gw\setminus \{0\}$, we derive
$$v(x)\leq \myfrac {C\gr(x) }{A^{1/q+1-p}\abs x^{(q+1)/(q+1-p)}},
$$
which implies (\ref{up-est2}).\qeda.\\

The construction of solutions of the quasilinear equations (\ref 
{equ+abs}) with
prescribed isolated singularity on the boundary of a general $C^2$ bounded domain 
$\Gw$ is settled upon similar constructions when the domain is either a half space, or a ball. 
%%%%%%%%%%%%%%%%%%%%%%%%%%%%%%%%%%%%%%%%%%%%%%%PROPOSITION HALF SPACE %%%%%%%%%%%%%%%%%%%%%%%%%%
%%%%%%%%%%%%%%%%%%%%%%%%%%%%%%%%%%%%%%%%%%%%%%%%%%%%%%%%%%%%%%%%%%%%%%%%%%%%%%%%%%%%%%
\bprop {halfspace} Assume $N-1<q<2N-1$ and let 
$H=\BBR^N_{+}=\{x=(x_{1},...,x_{N}):x_{N}>0)\}$ and  $k>0$. 
Then there exists a unique positive
solution $u=u^{{ H}}_{k}\in C^{1}(\overline {H}\setminus\{0\})$ of (\ref {equ+abs}) in 
$H$ which vanishes on 
$\prt H\setminus\{0\}$ and satisfies , 
\begin {equation}\label {sing-k}
u(x)=k\myfrac {x_{N}}{\abs x^{2}}(1+\circ (1)) \mbox { as }x\to 
0.
\end {equation}
\es
\Proof Since the function $x\mapsto kx_{N}\abs x^{-2}$ is $N$-harmonic in $H$ and 
vanishes on 
$\prt H\setminus\{0\}$, it is a supersolution of (\ref {equ+abs}). 
We write spherical coordinates in $\BBR^N$ under the form
\begin {equation}\label {spher}
x=\left\{(r,\gs)\in[0,\infty)\ti S^{N-1}=(r,\sin\gf\,\gs',\cos\gf):\gs'\in S^{N-2},
\gf\in [0,\gp]\right\},
\end {equation}
then 
$$Du=u_{r}{\bf i}+\myfrac {1}{r}\nabla_{\gs}u,
$$
where ${\bf i}=x/\abs x$, $\nabla_{\gs}$ denotes the covariant gradient on $S^{N-1}$, 
and equation (\ref {equ+abs}) 
takes the form 
\begin {equation}\label {psrad1}\BA{l}
-r^{1-N}\left(r^{N-1}\left(u_{r}^2+r^{-2}\abs {\nabla_{\gs}u}^2\right)^{(N-2)/2}
u_{r}\right)_{r}\\[4mm]
\qquad\qquad -r^{-2}div_{\gs}.\left(\left(u_{r}^2+
r^{-2}\abs {\nabla_{\gs}u}^2\right)^{(N-2)/2}\nabla_{\gs}u\right)
+\abs u^{q-1}u=0.
\EA\end {equation}
Next
$$\nabla_{\gs}u=-u_{\gf}{\bf e}+\myfrac {1}{\sin\gf}\nabla_{\gs'}u
$$
where ${\bf e}$ is derived from $x/\abs x$ by the rotation with angle $\gp/2$ 
in the plane $0,x,\rm N$ ($\rm N$ being the North pole), and $\nabla_{\gs'}$ is the covariant 
gradient on $S^{N-2}$ and (see \cite {BV})
\begin {equation}\label {psrad2}\BA{l}
div_{\gs}.\left(\left(u_{r}^2+
\myfrac{\abs {\nabla_{\gs}u}^2}{r^{2}}\right)^{(N-2)/2}\nabla_{\gs}u\right)\\
\phantom{-----}=\myfrac {1}{\sin^{N-2}\gf}\left(\sin^{N-2}\gf\left(u_{r}^2+\myfrac {u_{\gf}^2}{r^{2}}+\myfrac {\abs{\nabla_{\gs'}u}^2}{r^{2}\sin^2\gf}\right)^{(N-2)/2}\!\!u_{\gf}\right)_{\gf}
\\[6mm]
\phantom{-----}
+\myfrac {1}{\sin^2\gf}\,div_{\gs'}\left(\left(u_{r}^2+
\myfrac {u_{\gf}^2}{r^{2}}+
\myfrac {\abs{\nabla_{\gs'}u}^2}{r^{2}\sin^2\gf}\right)^{(N-2)/2}
\!\!\nabla_{\gs'}u\right).
\EA\end {equation}
 If $u$ depends only on $r$ and $\gf$, (\ref {equ+abs}) takes the form
\begin {equation}\label {psrad3}\BA{l}
-r^{1-N}\left(r^{N-1}\left(u_{r}^2+r^{-2}u_{\gf}^2\right)^{(N-2)/2}u_{r}\right)_{r}\\
\qquad\qquad-r^{-2}\sin^{2-N}\gf\left(\sin^{N-2}\gf\left(u_{r}^2+r^{-2}
u_{\gf}^2\right)^{(N-2)/2}u_{\gf}\right)_{\gf}
+\abs u^{q-1}u=0.
\EA\end {equation}
{\it Step 1} We look for a local subsolution $w$ under the form
$$w(r,\gs)=k(1-r^{\ga})r^{-1}\cos\gf\qquad r>0\,,\;\gf\in [0,\gp/2].
$$
where $\ga>0$ is to be determined. Then 
$$w_{r}=-kr^{-2}(1+(\ga-1)r^{\ga})\cos\gf\,\mbox { and }\;w_{\gf}
=-kr^{-1}(1-r^{\ga})\sin\gf$$
$$w_{r}^2+r^{-2}w_{\gf}^2:=P=k^2r^{-4}\left(1+2(\ga\cos^2\gf-1)
r^{\ga}+r^{2\ga}((\ga^2-2\ga)\cos^2\gf+1)\right)$$ 
$$w_{rr}=kr^{-3}(2-(\ga-1)(\ga-2)r^{\ga})\cos\gf\,\mbox { and }\;w_{\gf\gf}
=-kr^{-1}(1-r^{\ga})\cos\gf$$
$$P_{r}=-2k^2r^{-5}\left[2+(4-\ga)(\ga\cos^2\gf-1)r^{\ga}+
(2-\ga)((\ga^2-2\ga)\cos^2\gf+1)r^{2\ga}\right]
$$
$$P_{\gf}=-k^2\ga r^{\ga-4}\left[2+(\ga-2)r^\ga\right]\sin 2\gf,
$$
$$P_{r}w_{r}+r^{-2}P_{\gf}w_{\gf}=2k^3r^{-7}\left[2+(5\ga-6+
(2\ga-\ga^2)\cos^2\gf)r^{\ga}+O(r^{2\ga}))\right]\cos\gf,
$$
$$\BA {l}(N-1)r^{-1}w_{r}+w_{rr}+(N-2)r^{-2}\cot\gf\,w_{\gf}+r^{-2}w_{\gf\gf}\\
\phantom{(N-1)r^{-1}w_{rr}+(N-2)r^{-2}\cot\gf\,w_{\gf}+}=
kr^{-3}\left[4-2N)+(2-\ga)(N+\ga-2)r^\ga\right]\cos\gf.
\EA$$
Since
$$\BA {l}
-div\left(\abs {Dw}^{N-2}Dw\right)+w^q=Lw\\[2mm]
\phantom {--------}
=-P^{(N-2)/2}\left[(N-1)r^{-1}w_{r}+w_{rr}+(N-2)r^{-2}\cot\gf\,w_{\gf}+
r^{-2}w_{\gf\gf}\right]\\[2mm]
\phantom {-------------------}
-\myfrac {N-2}{2}P^{(N-4)/2}\left[P_{r}w_{r}+r^{-2}P_{\gf}w_{\gf}\right]+w^q,
\EA$$
and
$$w^q=k^q(1-r^{\ga})^qr^{-q}\cos^q\phi=k^q(1-qr^{\ga}+O(r^{2\ga}))r^{-q}\cos^q\gf,
$$
a straightforward computation leads to
\begin {equation}\label {comp}\BA {l}
Lw=k^{p-1}\ga\left[3-2N+(2+\ga)(N-2)\cos^2\gf+
O(r^{\ga})\right]P^{(N-4)/2}r^{\ga-7}\cos\gf\\[2mm]
\phantom {-------------------}+k^q(1-qr^{\ga}+O(r^{2\ga}))r^{-q}\cos^q\gf\\[2mm]
\phantom {Lw}
=k^{p-1}\ga\left[3-2N+(2+\ga)(N-2)\cos^2\gf\right]r^{-(2N-1)+
\ga}\cos\gf+k^qr^{-q}\cos^q\gf\\[2mm]
\phantom {Lw------}
-qk^qr^{-q+\ga}\cos^q\gf
+O(r^{-(2N-1)+2\ga}\cos\gf)+O(r^{-q+2\ga}\cos\gf).
\EA\end {equation}
By assumption $q<2N-1$. If we choose $\ga<\min\{2N-1-q,1/(N-2)\}$, 
there exists $R\in (0,1]$ such that $Lw\leq 0$ on $H\cap B_{R}$.\medskip

\nind {\it Step 2} Next we construct a solution $u_{R}$ in $B_{R}\cap H$  
which vanishes on
$\prt B_{R}\cap H$ and on $\prt H\setminus\{0\}$ and satisfies
\begin {equation}\label {asym1}
\lim_{r\to 0}\myfrac {ru_{R}(r,\gs)}{\cos\gf}=k.
\end {equation}
Let $\ell_{R}=k(1-R^\ga)R^{-1}$. Inasmuch  $w-\ell_{R}$ is a subsolution, for 
any $\ge>0$ we can construct a nonnegative solution $u_{\ge}$ of 
(\ref {equ+abs}) in $H\cap(B_{R}\setminus B_{\ge})$ which vanishes 
on $H\cap \prt B_{R}$ and on $\prt H \cap (B_{R}\setminus B_{\ge})$ 
and takes the value $k\ge^{-2}x_{N}$ on $H\cap\prt B_{\ge}$. By comparison
\begin {equation}\label {asym2}
(w(x)-\ell_{R})_{+}\leq u_{\ge}(x)\leq kx_{N}\abs x^{-2}.
\end {equation}
Furthermore, $\ge\mapsto u_{\ge}$ is increasing. Set $u=u_{R}=\lim_{\ge\to 0}u_{\ge}$, then $u$ is a solution of (\ref {equ+abs}) in $H\cap B_{R}$ which vanishes on $\prt B_{R}\cap H$ and on $\prt H\setminus\{0\}$ and satisfies the same inequality (\ref {asym2}) as  $u_{\ge}$, but in whole
 $H\cap B_{R}$. This implies that (\ref {asym1}) holds uniformly on $[0,\gp/2-\gd]$, for any $\gd>0$. In order to improve this inequality, we perform a scaling: for $r>0$, we set $u^r(x)=ru(rx)$. Then $u^r$ satisfies
 \begin {equation}\label {scale}
 -div\left(\abs {Du^r}^{N-2}Du^r\right)+r^{2N-1-q}(u^r)^q=0
\end {equation}
 in $H\cap B_{R/r}$ where there holds
\begin {equation}\label {asym3}
k(x_{N}\abs x^{-2}(1-r^\ga\abs x^{\ga})-\ell_{R})_{+}\leq u^r(x)\leq kx_{N}\abs x^{-2}.
\end {equation}
Since $u^r$ is uniformly bounded for $1/2\leq\abs x\leq 2$, it follows from regularity theory  \cite {Lie} that it is also bounded in the $C^{1,\ga}$-topology of $2/3\leq\abs x\leq 3/2$. Using Ascoli's theorem and the fact that $u^r(x)$ converges to $kx_{N}\abs x^{-2}$ pointwise and locally uniformly, it follows that $Du^r(x)=r^2Du(rx)$ converges uniformly in $\{x\in H:2/3\leq\abs x\leq 3/2\}$ to 
$-2kx_{N}\abs x^{-4}x+k\abs x^{-2}{\bf e}_N$ which is the gradient of $x\mapsto kx_{n}\abs x^{-2}$. Using the expression of $Du$ in spherical coordinates we obtain
$$r^2u_{r}{\bf i}-ru_{\gf}{\bf e}+\myfrac {r}{\sin\gf}\nabla_{\gs'}u\to
-2k\gs_{N}{\bf i}+k{\bf e}_N\mbox { uniformly on }S^{N-1}_{+}\mbox { as }r\to 0,
$$
where $\gs_{N}=\langle\gs,{\bf e}_N\rangle$. Inasmuch ${\bf i}$, ${\bf e}$ and $\nabla_{\gs'}u$ are orthogonal, the component of ${\bf e}_N$ is $\sin\gf$, thus
\begin {equation}\label{rep1}
ru_{\gf}(r,\gs',\gf)\to -k\sin\gf\;\mbox { as }r\to 0.
\end {equation} 
Since
\begin {equation}\label{rep2}
u(r,\gs',\gf)=\myint{\gp/2}{\gf}u_{\gf}(r,\gs',\gth)\,d\gth,
\end {equation}
the previous convergence estimate establishes (\ref {asym1}).\medskip

\nind {\it Step 3} Construction of the solution in $H$. Let $\eta$ be the truncation function introduced in the proof of \rprop {KOV}, and $\eta_{\ge}(r)=\eta((r-\ge)_{+})$. Then the function 
$u_{R,\ge}$ defined by $u_{R,\ge}=\eta_{\ge}\circ u_{R}$ in $H\cap B_{R}$ and zero outside, is a subsolution of (\ref {equ+abs}) in $H$ which vanishes on $\prt H\setminus \{0\}$ and satisfies (\ref {asym1}). Using the same device as in Step 2, we construct a sequence of solutions $u_{\gd}$ ($\gd>0$) of (\ref {equ+abs}) in $H\setminus B_{\gd}$ with boundary value $k\gd^{-2}x_{N}$ on $\prt B_{\gd}\cap H$, zero on $\prt H\setminus B_{\gd}$ and satisfies
$$u_{R,\ge}\leq u_{\gd}\leq kx_{N}\abs x^{-2}.
$$
When $\gd\to 0$, $u_{\gd}$ decreases and converges to some $u$ which satisfies (\ref {equ+abs}) and  the previous inequality. Letting successively $\ge\to 0$ and $\eta(r)\to r_{+}$ we obtain that $u$ satisfies
\begin {equation}\label {asym4}
\check u_{R}(x)\leq u(x)\leq kx_{N}\abs x^{-2}\quad\mbox {in }H,
\end {equation}
where $\check u$ is the extension of $u$ by zero outside $B_{R}$. The proof of (\ref {sing-k}) is the same as in Step 2.\medskip

\nind {\it Step 4} Uniqueness. Let $u$ and $\hat u$ be two solutions of (\ref{equ+abs}) satisfying (\ref{sing-k}) and $\ge>0$. Then $u_{\ge}=(1+\ge)u+\ge$ is a super solution which is positive of $\prt H\setminus\{0\}$. Inasmuch it dominates $\hat u$ both in a neighborhood of $0$ and in a neighborhood of infinity, it dominates $\hat u$ in $H$. Letting $\ge\to 0$ yields to $u\geq \hat u$. Similarly $\hat u\geq u$.\qeda
%%%%PROPOSITION BALL %%%%%%%%%%%%%%%%%%%%%%%%%%
 %%%%%%%%%%%%%%%%%%%%%%%%%%%%%%%%%%%%%%%%%%%%%%%
%%%%%%%%%%%%%%%%%%%%%%%%%%%%%%%%%%%%%%%%%%%%%%%
\bprop {ball}Assume $N-1<q<2N-1$ and let $B=B_{1}(0)$, $a\in \prt B$ and $k>0$. Then there exists a unique 
function $u=u^B_{k,a}\in C^{1}(\overline {B}\setminus\{a\})$ which vanishes on 
$\prt B\setminus\{a\}$ and satisfies (\ref {equ+abs}) in 
$B$ and
\begin {equation}\label {sing-k2}
u(x)=k\myfrac {1-\abs x}{\abs {x-a}^{2}}(1+\circ (1)) \mbox { as }x\to 
a.
\end {equation}
\es
\Proof With a change of coordinates, we can assume that $B$ has center $m=(0,...,0,-1/2)$ and $a$ is the origin of coordinates. We denote by $\gw$ the point $(0,...,0,-1)$ and by $\CI_{\gw}$ the inversion with center $\gw$ and power $1$. By this involutive transformation, the half space $H=\{x\in\BBR^N :x_{N}>0\}$ is transformed into the ball $B^*=\{x\in\BBR^N:\abs x^2+x_{N}<0\}$. Thus the function $x\mapsto P_{k}(x)= -k(\abs x^2+x_{N})/2\abs x^2$ is $N$-harmonic and positive  in $B^*$, vanishes on $\prt B^*\setminus \{0\}$ and is singular at $0$. 
Let $v_{k}$ be the solution of (\ref{equ+abs}) in $H$ satisfying (\ref{sing-k}), and $u_{k}=v_{k}\circ \CI_{\gw}$. Then $u_{k}\in C (\overline{B^*}\setminus\{0\})$ satisfies
\begin {equation}\label {equ+abs2}\left\{\BA{l}
-div\left(\abs {Du_{k}}^{N-2}Du_{k}\right)+\abs {x-\gw}^{-2N}u^q_{k}=0\quad\mbox {in } 
B^*\\
\phantom {-div\left(\abs {Du_{k}}^{N-2}Du_{k}\right)+\abs {x-\gw}^{-2N}}
u_{k}=0\quad\mbox {on } \prt B^*\setminus\{0\}.
\EA\right.\end {equation}
Furthermore $u_{k}\leq P_{k}$ and
\begin {equation}\label {asym5}
P_{k}(x)=k\myfrac {1/4-\abs{x-m}^2}{2\abs x^2}=k\myfrac{1/2-\abs {x-m}}{2\abs x^2}
(1+\circ (1))=u_{k}(x)(1+\circ (1))
\end {equation}
as $x\to 0$. Inasmuch $\abs {x-\gw}\leq 1$, $u_{k}$ is a subsolution of (\ref {equ+abs}) in $B^*$.
For $\ge>0$ we construct a solution $v_\ge$ of (\ref {equ+abs}) in  $B^*\setminus B_\ge(0)$ with boundary value $P_k$. By the maximum principle $u_{k}\leq v_\ge\leq P_{k}$ in $B^*\setminus B_\ge(0)$. Since the sequence $\{v_\ge\}$ is monotone,  
we obtain that there exists a solution $\lim_{\ge\to 0} v_\ge:= u\in  C^1 (\overline{B^*}\setminus\{0\})$ of (\ref {equ+abs}) 
in $B^*$ which satisfies 
\begin {equation}\label {asym6}
u_{k}(x)\leq u(x)\leq P_{k}(x)\quad\mbox {in }B^*,
\end {equation}
and
\begin {equation}\label {asym7}
u(x)=k\myfrac{1/2-\abs {x-m}}{2\abs x^2}(1+\circ (1)).
\end {equation}
We change the variables in setting $x'_{N}=x_{N}+1/2$ and $x'_{i}=x_{i}$ ($i=1,...,N-1$). 
We define $u'(x')=u(x)$ and denote by $a$ the point $(0,...,0,1)$. Clearly
$u'$ 
satisfies (\ref {equ+abs}) in $B_{1/2}$, vanishes on $\prt B_{1/2}\setminus\{a\}$ and
\begin {equation}\label {asym7'}
u'(x)=k\myfrac {1/2-\abs x}{2\abs {x-a/2}^{2}}(1+\circ (1)) \mbox { as }x\to 
a/2.
\end {equation}
By the transformation $\ell\mapsto \ell^{p/(q+1-p)}u'_{k}(\ell x)$, where $\ell=1/2$, we obtain a 
solution $u_{k,a}$ of (\ref {equ+abs})  in $B$ which verifies
\begin {equation}\label {asym8}
u_{k,a}(x)=2^{N/(q+1-N)}k\myfrac {1-\abs x}{\abs {x-a}^{2}}(1+\circ (1)) \mbox { as }x\to 
a.
\end {equation}
Because $k$ is arbitrary, (\ref {sing-k2}) follows. Uniqueness of the solution is obtained as in \rprop {halfspace} with $u_\ge=(1+\ge)u$.\qeda 
%%%%%%%%%%%%%%%%%%%%%%%%%%%%%%%%%%%%%%%%%%%%%%%
%%%%PROPOSITION COMPLEMENT OF BALL %%%%%%%%%%%%%%%%%%%%%%%%%%
%%%%%%%%%%%%%%%%%%%%%%%%%%%%%%%%%%%%%%%%%%%%%%%
\bprop {compball} Assume $N-1<q<2N-1$ and let $G=\overline{B}^c$, $a\in \prt B$ and $k>0$. Then there exists a unique 
function $u=u^{B^c}_{k,a}\in C^{1}(\overline {G}\setminus\{a\})$ which vanishes on 
$\prt B\setminus\{a\}$ and satisfies (\ref {equ+abs}) in 
$G$ and
\begin {equation}\label {sing-k'2}
u(x)=k\myfrac {\abs x-1}{\abs {x-a}^{2}}(1+\circ (1)) \mbox { as }x\to 
a.
\end {equation}
\es
\Proof Uniqueness follows from (\ref{sing-k'2}) by the same method as in \rprop {halfspace} and \rprop {ball}. Actually, it will be proved in \rth{gendom}.  For existence we perform the inversion $\CI^{1}_{0}$ with center $0$ and 
power $1$. It transforms the function $u^{B}_{k,a}$ constructed in the 
previous proposition into a function $v\in C^{1}(\overline {G}\setminus\{a\})$ which vanishes on 
$\prt B\setminus\{a\}$ and satisfies (\ref{sing-k'2}). Furthermore $v$ 
is solution of 
\begin {equation}\label {equ+abs'}
-div\left(\abs {Dv}^{N-2}Dv\right)+\abs x^{-2N}\abs v^{q-1}v=0
\end {equation}
in $G$. Since $\abs x>1$, $v$ is a super solution for (\ref {equ+abs}) 
in $G$. With no loss of generality we can assume that $a=(0,...0,1)$ and let $u_{k,a}^{H^1}$ be the solution of (\ref {equ+abs}) in $H^1=\{x=(x_1,...,x_N:x_N>1)\}$ 
satisfying (\ref{sing-k}) already 
constructed in \rprop{halfspace}. Then $\gu_\ge=\eta(u_{k,a}^{H^1})$ is a subsolution in 
$G$ (where $\eta_\ge$ has been defined in the proof of \rprop {KOV}). By the same approximation as in 
the previous proposition, we construct an increasing sequence 
$\{u_{\ge}\}$ ($\ge>0$) of 
solutions of (\ref{equ+abs}) in $G\setminus B_{\ge}(a)$ which vanishes 
on $\prt G\setminus B_{\ge}(a)$, takes the value $v$ on $ G\cap \prt 
B_{\ge}(a)$ and verifies $\gu_\ge\leq u_{\ge}\leq v$ in $G\setminus B_{\ge}(a)$. Letting $\ge\to 0$, we obtain the existence of a solution $ u^*$ in $G$ which satisfies 
\begin {equation}\label {comp'2}
\tilde u_{k,a}^{H^1}\leq u^*\leq v\quad\mbox {in }G
\end {equation}
where we denote by $\tilde u_{k,a}^{H^1}$ the extension of $ u_{k,a}^{H^1}$ by zero in $\overline {H^1}^c$. We conclude that (\ref {sing-k'2}) holds in $H^1$. In order to extend this convergence to whole $G$, we proceed as in the proof of \rprop{halfspace}, with a minor modification due to the geometry. We put the origin of coordinates at $a$, takes the same spherical coordinates and obtain again that
$$r^2u^*_r{\bf i}-ru^*_\gf{\bf e}+\myfrac{r}{\sin\phi}\,\nabla_{\gs'}u^*\to -2k\gs_N{\bf i}+k{\bf e}\quad\mbox {  uniformly on $S^{N-1}_+$ as }r\to 0.
$$
Therefore (\ref{rep1}) holds for any $\gf\in [0,\gp/2]$. For $r>0$, the angle $\gf$ ranges from $\psi(r)=\cos^{-1}(-r/2)$ to $0$
(here is the difference with the half-space case) and $\abs x^2\nabla u(x)$ remains bounded in this domain, by the regularity theory for quasilinear elliptic equations. Since
 \begin {equation}\label{rep2'}
u^*(r,\gs',\gf)=\myint{\psi(r)}{\gf}u^*_{\gf}(r,\gs',\gth)\,d\gth,
\end {equation}
we derive, as in the proof of \rprop {halfspace}, 
 \begin {equation}\label{rep3}
\lim_{r\to 0}u^*(r,\gs',\gf)=k\cos\gf\quad\mbox{ uniformly on } [0,\gp/2].
\end {equation}
The proof that (\ref{sing-k'2}) holds is a particular case of \rth {gendom}.\qeda \\

%%%%%%%%%%EXTENSION LEMMA%%%%%%%%%%%%%%%%%%%%%%%%%%%%%%
%%%%%%%%%%%%%%%%%%%%%%%%%%%%%%%%%%%%%%%%%%%%%%%%
%%%%%%%%%%%%%%%%%%%%%%%%%%%%%%%%%%%%%%%%%%%%%%%%%
 In a general domain we have to extend the solution through the boundary. We denote by $\dot\gr(x)$ the signed distance from $x\to\prt\Gw$, that is $\dot\gr(x)=\gr(x)$ if $x\in\Gw$ and $\dot\gr(x)=-\gr(x)$ if $x\in \Gw^c$. Since $\prt\Gw$ is $C^2$, there exists $\gb_{0}>0$ such that if $x\in \BBR^N$ verifies 
 $-\gb_{0}\leq\dot\gr(x)\leq\gb_{0}$, there exists a unique $\xi_{x}\in\prt\Gw$ such that 
 $\abs{x-\xi_{x}}=\abs{\dot\gr(x)}$. Furthermore, if $\gn_{\xi_{x}}$ is the outward unit vector to $\prt\Gw$ at $\xi_{x}$, $x=\xi_{x}-\dot\gr(x)\gn_{\xi_{x}}$. In particular $\xi_{x}-\dot\gr(x)\gn_{\xi_{x}}$ and $\xi_{x}+\dot\gr(x)\gn_{\xi_{x}}$ have the same orthogonal projection $\xi_{x}$ onto $\prt\Gw$.\smallskip
 
 Let $T_{\gb_{0}}(\Gw)=\{x\in\BBR^N:-\gb_{0}\leq\dot\gr(x)\leq\gb_{0}\}$, then the mapping 
 $\Gp:[-\gb_{0},\gb_{0}]\ti\prt\Gw\mapsto T_{\gb_{0}}(\Gw)$ defined by $\Gp(\gr,\xi)=\xi-\dot\gr{\bf\gn}(\xi)$ is a $C^2$ diffeomorphism. Moreover $D\Gp(0,\xi)(1,e)=e-\gn_{\xi}$ for any $e$ belonging to the tangent space $T_{\xi}(\prt\Gw)$ to $\prt\Gw$ at $\xi$. If $x\in T_{\gb_{0}}(\Gw)$, we define the reflection of $x$ through $\prt\Gw$ by
 $\psi(x)=\xi_{x}+\dot\gr(x)\gn_{xi_{x}}$. Clearly $\psi$ is an involutive diffeomorphism from $\overline\Gw\cap T_{\gb_{0}}(\Gw)$ to $\Gw^c\cap T_{\gb_{0}}(\Gw)$. Furthermore for any $\xi\in\prt\Gw$, $D\psi (\xi)=\CS_{T_\xi(\prt\Gw)}$ is the symmetry with respect to the tangent space $T_{\xi}(\prt\Gw)$ to $\prt\Gw$ at $\xi$.
 If a function $v$ is defined in $\Gw\cap T_{\gb_{0}}(\Gw)$, we define $\tilde v$ in $\Gw^c\cap T_{\gb_{0}}(\Gw)$ by
 \begin {equation}\label {ext}
\tilde v(x)=\left\{\BA{l}v(x)\quad\qquad\mbox {if }x\in \Gw\cap T_{\gb_{0}}(\Gw)\\[2mm]
-v\circ\psi(x)\quad\mbox {if }x\in \Gw^c\cap T_{\gb_{0}}(\Gw).
\EA\right.\end {equation}
\bprop{ext} Let $v\in C^{1,\ga}(\overline\Gw\cap T_{\gb_{0}}(\Gw)\setminus\{0\})$ be a solution of (\ref {equ+abs}) in $\Gw\cap T_{\gb_{0}}(\Gw)$ vanishing on $\prt\Gw\setminus\{0\}$. Then 
$\tilde v\in C^{1,\ga}(T_{\gb_{0}}(\Gw)\setminus\{0\})$ is solution of a quasilinear equation
 \begin {equation}\label {ext-equ}
-\mysum{j}{}\myfrac {\prt }{\prt x_{j}}\tilde A_{j}(x,D\tilde v)+\tilde b(x)\abs {\tilde v}^{q-1}\tilde v=0
\end {equation}
in $T_{\gb_{0}}(\Gw)\setminus\{0\}$ where the $\tilde A_{j}$ and $\tilde b$ are $C^1$ functions defined in $T_{\gb_{0}}(\Gw)$ where they verify 
 \begin {equation}\label {ext1}\left\{\BA {l}
 (i)\quad \tilde A_{j}(x,0)=0\\[2mm]
 
 (ii) \quad\mysum{i,j}{}\myfrac {\prt \tilde A_{j}}{\prt \eta_{i}}(x,\eta)\xi_{i}\xi_{j}
 \geq \Gg\abs\eta^{p-2}\abs \xi^2
 \\
 
 (iii) \quad\mysum{i,j}{}\abs {\myfrac {\prt  \tilde A_{j}}{\prt \eta_{i}}(x,\eta)}\leq \Gg\abs\eta^{p-2}\\
 
 (iv)\quad \Gg\geq \tilde b(x)\geq \gg
\EA\right.
\end {equation}
for all $x\in T_{\gb}(\Gw)\setminus\{0\}$ for some $\gb\in (0,\gb_{0}]$, $\eta\in \BBR^N$, $\xi\in\BBR^N$ and some 
$0<\gg\leq \Gg$.
\es
\Proof The assumptions (\ref {ext1}) implies that weak solutions of (\ref {ext-equ}) are $C^{1,\ga}$, for some $\ga>0$ \cite {To1} and satisfy the standard a priori estimates. As it is defined the function $\tilde v$ is clearly $C^1$ in $T_{\gb_{0}}(\Gw)\setminus\{0\}$. Writing
$Dv(x)=-D(\tilde v\circ \psi(x))=-D\psi(x)(D\tilde v(\psi(x)))$ and $\tilde x=\psi(x)=\psi^{-1}(x)$
$$\BA {l}
\myint {\Gw\cap T_{\gb_{0}}(\Gw)}{}\left(\abs {Dv}^{p-2}Dv.D\gz +\abs v^{q-1}v\gz\right)dx\\

\phantom {-----}
=\myint {\overline\Gw^c\cap T_{\gb_{0}}(\Gw)}{}\left(\abs {D\psi(D\tilde v)}^{p-2}
D\psi (D\tilde v).D\psi(D\gz) +\abs {\tilde v}^{q-1}\tilde v\gz(\psi(\tilde x))\right)\abs {D\psi}d\tilde x.
\EA$$
But 
$$\BA {l}D\psi (D\tilde v).D\psi(D\gz)=\mysum{k}{}
\left(\mysum{i}{}\myfrac {\prt \psi_{i}}{\prt x_{k}}\myfrac {\prt \tilde v}{\prt x_{i}}\right)
\left(\mysum{j}{}\myfrac {\prt \psi_{j}}{\prt x_{k}}\myfrac {\prt \gz}{\prt x_{j}}\right)\\
\phantom{D\psi (D\tilde v).D\psi(D\gz)}
=\mysum{j}{}\left(\mysum{i,k}{}\myfrac {\prt \psi_{i}}{\prt x_{k}}\myfrac {\prt \psi_{j}}{\prt x_{k}}\myfrac {\prt \tilde v}{\prt x_{i}}\right)\myfrac {\prt \gz}{\prt x_{j}}.
\EA$$
We set $b(x)=\abs{D\psi}$, 
\begin {equation}\label {ext2}
A_j(x,\eta)=\abs {D\psi}\abs {D\psi(\eta)}^{p-2}\mysum{i}{}\left(\mysum{k}{}\myfrac {\prt \psi_{i}}{\prt x_{k}}\myfrac {\prt \psi_{j}}{\prt x_{k}}\right)\eta_i,
\end {equation}
and
\begin {equation}\label {ext3}
A(x,\eta)=(A_1(x,\eta),...,A_N(x,\eta))=\abs {D\psi}\abs {D\psi(\eta)}^{p-2}(D\psi)^tD\psi(\eta).
\end {equation}
For any $\xi\in\prt\Gw$, the mapping
$D\psi_{\prt\Gw}(\xi)$ is the symmetry with respect to the hyperplane $T_{\xi}(\prt\Gw)$ tangent to $\prt\Gw$ at $\xi$, so $\abs {D\psi(\xi)}=1$. Inasmuch $D\psi$ is continuous, a lengthy but standard computation leads to the existence of some $\gb\in (0,\gb_{0}]$ such that (\ref {ext1}) holds in $T_{\gb}(\Gw)\cap\overline\Gw^c$. If we define $\tilde A$ (resp. $\tilde b$) to be $\abs\eta^{p-2}\eta$ (resp $1$) on $T_{\gb}(\Gw)\cap\overline\Gw$ and $A$ (resp. $\abs{D\psi}$) on $T_{\gb}(\Gw)\cap\overline\Gw^c$, then inequalities (\ref {ext1}) are satisfied in $T_{\gb}(\Gw)$.\qeda \\

\noindent \Remark Notice that, similarly to the $p$-laplacian, the vector field $\tilde A$ is positively homogeneous with exponent $p-1$ with respect to $\eta$. Furthermore, if for  $r>0$ we set
$\tilde A^r_{j}(x,\eta)=\tilde A_{j}(rx,\eta)$ , then $\tilde A^r_{j}$ satisfies the same estimates (\ref {ext1}) as $A_{j}$, uniformly in 
$T_{r^{-1}\gb}(r^{-1}\Gw)$, for $0<r\leq 1$. Furthermore
$$\lim_{r\to 0}A^r_{j}(x,\eta)=\abs\eta^{p-2}\eta_j\;\forevery\eta\in\BBR^N\,,\;\;\forall j=1,...,N,
$$
and this limit is uniform on the bounded subsets of $\BBR^N$.
%%%%%%%%%%BIG THEOREM I%%%%%%%%%%%%%%%%%%%%%%%%%%%%%%
%%%%%%%%%%%%%%%%%%%%%%%%%%%%%%%%%%%%%%%%%%%%%%%%
%%%%%%%%%%%%%%%%%%%%%%%%%%%%%%%%%%%%%%%%%%%%%%%%%
\bth {gendom} Let $\Gw$ be a bounded domain with a $C^2$ boundary and $a\in\prt\Gw$. 
Assume $N-1<q<2N-1$ and denote by $\gr(x)$ the distance from $x$ to $\prt\Gw$. 
Then for any $k> 0$ there exists a unique function $u=u_{k,a}\in C(\overline \Gw\setminus \{a\})$ 
which vanishes on $\prt\Gw\setminus\{a\}$, is solution of (\ref {equ+abs}) and satisfies
\begin {equation}\label {sing-k3}
u_{k,a}(x)=k\myfrac {\gr(x)}{\abs {x-a}^{2}}(1+\circ (1)) \mbox { as }x\to a.
\end {equation}
\es
\Proof Uniqueness follows from (\ref {sing-k3}) by the same technique as in the previous 
propositions. For existence let $B^i_R$ be a ball of radius $R$ such that 
$B^i_R\subset\Gw$ and $a\in\prt B^i_R$, and let $\gw_i$ be its center. We denote 
by $U^i$ the solution of (\ref {equ+abs}) in $B^i_R$, which vanishes on 
$\prt B^i_R\setminus\{a\}$ and satisfies
\begin {equation}\label {sing-k4}
U^i(x)=k\myfrac{R-\abs {x-\gw_i}}{\abs{x-a}^2}(1+\circ (1))\mbox { as }x\to a.
\end {equation}
If we set $U_\gd=\eta_\gd(U^i)$, we have already seen that $\check U_\gd$, the extension 
of $U_\gd$ by zero outside its support, is a subsolution of (\ref {equ+abs}) in $\Gw$. 
Because $V^\Gw_{a}$, the $N$-harmonic function element of $C(\overline\Gw\setminus\{a\})$  vanishing on $\prt\Gw\setminus\{a\}$,  satisfies
\begin {equation}\label {sing-k5}
V^\Gw_{a}(x)=\myfrac {\gr(x)}{\abs {x-a}^{2}}(1+\circ (1)) \mbox { as }x\to a, \;x\in B^i_R,
\end {equation}
there holds $kV^\Gw_{a}\geq \check U_\gd$. If $\Gw_{\ge}=\Gw\setminus\{B_{\ge}(a)\}$ ($\ge>0$), 
we construct a solution $u_{\ge}\in C(\overline{\Gw_{\ge}})$ of (\ref {equ+abs}) in $\Gw_{\ge}$, 
which vanishes on $\prt\Gw\setminus B_{\ge}(a)$ and takes the value $kV^\Gw_{a}$ on  
$\prt B_{\ge}(a)\cap\Gw$. By the maximum principle $\ge\mapsto u_{\ge}$ is increasing and
$\check U_\gd\leq u_{\ge}\leq kV^\Gw_{a}$ in $\Gw_{\ge}$. Letting $\ge\to 0$ we obtain that 
$u_{\ge}$ converges in the $C_{loc}^{1}$-topology of $\overline {\Gw}\setminus\{a\}$ to 
a solution $u=u_{k,a}$ of (\ref{equ+abs}) in $\Gw$. It follows from the previous inequalities that
 \begin {equation}\label {sing-k6}
\check U_\gd(x)\leq u(x) \leq kV^\Gw_{a}(x) \forevery x\in \overline {\Gw}\setminus\{a\}.
\end {equation}
In order to prove the asymptotic behaviour, we proceed as in \rprop {ball} with the help of the reflection principle of \rprop{ext}. We fix the origin of 
coordinates at $a=0$ and the normal outward unit vector at $a$ to be $-{\bf e}_{N}$. If $\tilde u$ is the extension of $u$ by reflection through $\prt\Gw$, it satisfies
 \begin {equation}\label {ext-k1}
-\mysum {j}{}\myfrac {\prt }{\prt x_{j}}\tilde A_{j}(x,D\tilde u)+\tilde b(x)\abs {\tilde u}^{q-1}\tilde u=0
\end {equation}
in $T^\gb(\Gw)\setminus\{0\}$. For $r>0$, set $\tilde u^r(x)=r\tilde u(rx)$. Then $\tilde u^{r}$ is solution of 
  \begin {equation}\label {ext-k2}
-\mysum {j}{}\myfrac {\prt }{\prt x_{j}}\tilde A^r_{j}(x,D\tilde u^r)+r^{2N-1-q}\tilde b(rx)\abs {\tilde u^r}^{q-1}\tilde u^r=0
\end {equation}
in $T^{\gb r^{-1}}(\Gw^r)\setminus\{0\}$, where $\Gw^r:=r^{-1}\Gw$. By \cite [Th 2.4]{BV} there exists $C>0$ such that
$$kV^\Gw_{0}(x)\leq Ck\myfrac {\gr(x)}{\abs {x}^{2}}.
$$
Furthermore, for any $x\in T^{\gb}(\Gw)\setminus\{0\}$, $\gr(x):=\dist (x,\Gw)=\gr(\psi(x))$ (we recall that $\psi(x)$ is the symmetric of $x$ with respect to $\prt\Gw$ as it is defined in \rprop{ext}), and $c\abs x\leq \abs{\psi(x)}\leq c^{-1}\abs x$ for some $c>0$, the same relations holds if $T^{\gb}(\Gw)$ is replaced by $T^{\gb r^{-1}}(\Gw^r)$ and $\gr(x)$ by $\gr_{r}(x):=\dist (x,\Gw^r)$. Since $\Gw$ is $C^2$,
 $$\lim_{r\to 0}\myfrac {\gr(rx)}{r\gr_{r}(x)}=1$$
uniformly on bounded subsets of $\BBR^N$. Consequently
$$\abs{\tilde u^r}(x)\leq Ckr^{-1}\myfrac {\gr(rx)}{\abs {x}^{2}}= Ck\myfrac {\gr_{r}(x)}{\abs {x}^{2}}
(1+\circ(1)).
$$
For $0<a<b$ fixed and $r\leq r_{0}$ (for some $r_{0}\in (0,1]$) the spherical shall $\Gg_{a,b}=\{x\in \BBR^N:a\leq\abs x\leq b\}$ is included into $ T^{\gb r^{-1}}(\Gw^r)$. By the classical regularity theory for quasilinear equations \cite {To1} and \rprop{ext},  there holds
 \begin {equation}\label {sing-k7}
\norm {D\tilde u^r}_{C^\ga (\Gg_{2/3,3/2})}\leq C_{r}\norm {\tilde u^r}_{L^\infty (\Gg_{1/2,2})},
\end {equation}
where $C_{r}$  remains 
bounded because $r\leq 1$. By Ascoli's theorem and  (\ref{sing-k6}) $\tilde u^r(x)$ converges to 
$kx_{N}\abs x^{-2}$ in the $C^1 (B_{3/2}\setminus B_{1/2})$-topology. This implies in 
particular
$$\lim_{r\to 0}r^2D\tilde u(rx)=-2kx_{N}x\abs{x}^{-4} + k \abs{x}^{-2}{\bf e}_{N}. 
$$
If we take in particular $\abs x=1$, we derive
\begin {equation}\label {asym'1}
\lim_{r\to 0}(r\tilde u(r,\gs),r^2\nabla \tilde u(r,\gs))=(k\cos\gf,-k\sin\gf \bf{e}_N ),
\end {equation}
uniformly with respect to $\gs=(\sin\gf\,\gs',\cos\gf)\in S^{N-2}\ti [0,\gp]$. 
Because $\prt\Gw$ is $C^2$ there exists $\ge_{0}>0$ and a $C^2$ real valued function $h$ defined in 
$\Gth_{\ge_{0}}:=B_{\ge_0}\cap \prt H$ (we recall that $\prt H=\{x=(x',0)\}$) and an open neighborhood $\CV_{\ge_{0}}$ of $0$ such that 
$\prt\Gw\cap \CV_{\ge_{0}}=\{x=(x',x_{N}:x_{N}=h(x')\}$, and $Dh(0)=0$ (this expresses the fact that $\prt H=T_{0}(\prt\Gw)$). If we define $\Psi$ by 
$$\Psi(x)=(x',x_{N}-h(x'))\forevery x\in \CV_{\ge_{0}}.
$$
then $det (D\Psi)=1$ and $D\Psi(0)=I$. Up to replacing $\ge_{0}$ by a smaller quantity, $\Psi$ is a $C^2$ diffeomorphism from $\CV_{\ge_{0}}$ into a neighborhood $\CV'$ of $0$ such that $\CV_{\ge_{0}}\cap\prt\Gw)=\Gth_{\ge_{0}}$. Because $\dist (\Psi(x),\prt H)=x_{N}-h(x')$, $\dist (\Psi(x),\prt H)=\gr(x)(1+\circ (1))$ as $x\to 0$. Thus, if
we set $x=\Psi^{-1}(y)$ and $\tilde u(x)=u^*(y)$, (\ref{asym'1}) is equivalent to
\begin {equation}\label {asym'2}
\lim_{\abs y\to 0} (\abs y u^*(\abs y,\gs),\abs y^2\nabla u^*(\abs y,\gs))=(k\cos\gf,-k\sin\gf \,\bf{e}_N ),
\end {equation}
uniformly on $S^{N-1}$, thus
\begin {equation}\label {asym'3}
\abs y u^*(\abs y,\gs)=k\sin\gf \,(1+\circ (1))\quad \mbox {as }\,\abs y\to 0
\end {equation}
uniformly with respect to $\gs\in S^{N-1}_+$, because $u^*$ vanishes on $B_{\ge_0}\cap\prt H\setminus\{0\}$. This implies (\ref{sing-k3}).
\qeda\\
%%%%%%%%%%%%%%%%%%%%%%%%%%%%%%%%%%%%%%%%%%%%%%%%%%%%%%%%%%%%TRANSITION TO STRONG SING.%%%%%%%%%%%%%%%%%%%%%%%%%%%

Clearly the mapping $k\mapsto u_{k,a}$ is increasing. As $u_k$ satisfies the estimates 
(\ref {up-est2}) and (\ref{sing-k'2}), $u_{k,a}$ converges in the $C^1_{loc}(\overline\Gw\setminus\{a\})$-topology, 
as $k\to\infty$, to some $u_{\infty,a}$, solution of (\ref {equ+abs}) in $\Gw$, 
vanishes on $\prt\Gw\setminus\{a\}$ and satisfies
\begin {equation}\label {sing-inf}
\lim_{x\to a} \myfrac {\abs {x-a}^{2}u_{\infty,a}(x)}{\gr(x)}=\infty.
\end {equation}
In order to describe the precise behaviour of $u_{\infty,a}$, we have to introduce separable 
solutions of (\ref {equ+abs}) in $\BBR^N\setminus\{0\}$: if we look for solutions $u$ 
under the form $u(r,\gs)=r^\gb\gw(\gs)$, then $\gb=-\gb_{q}=-N/(q+1-N)$ and $\gw$ satisfies
\begin {equation}\label {sep1}
-div_{\gs}\left(\left(\gb^2_{q}\gw^2+\abs {\nabla_{\gs}\gw}^{2}\right)^{(N-2)/2}
\nabla_{\gs}\gw\right)-\Gl\left(\gb^2_{q}\gw^2+\abs {\nabla_{\gs}\gw}^{2}
\right)^{(N-2)/2}\gw+\abs \gw^{q-1}\gw
=0
\end {equation}
on $S^{N-1}$ where $\Gl=(N-1)\gb_{q}^2$. We shall denote by $\CS_{q}$  the set of
(always $C^{1,\ga}$) solutions of (\ref {sep1}). If $u$ is a separable solution of (\ref {equ+abs}) 
in $H$ which vanishes on $\prt H\setminus\{0\}$, the function $\gw$ is a solution of (\ref {sep1}) 
in $S^{N-1}_{+}$ which vanishes on $\prt S^{N-1}_{+}=S^{N-2}$. We shall denote by 
$\CS^*_{q}$ the set of such functions and by $\CS^*_{q\,+}$ the subset of positive solutions. 
We recall some simple facts
%%%%%%%%%%%%%%%%%%%%%%%%%%%%%%%%%%%%%%%%%%%%%%%%
%%%%%%%%%%STRONG SINGULARITY%%%%%%%%%%%%%%%
%%%%%%%%%%%%%%%%%%%%%%%%%%%%%%%%%%%%%%%%%%%%%%%%%%%
\bprop{psirad} (i) For any $q>N-1$, $\CS_{q}$ contains at least the three constant functions $0$ 
and
$\pm ((N-1)\gb_q^N)^{1/(q+1-N)}$.\smallskip

\nind (ii) For any $q\geq 2N-1$, $\CS^*_{q}=\{0\}$.\smallskip

\nind (iii) For any $q\in (N-1,2N-1)$, $\CS^*_{q\,+}$ contains a unique element.
\es
\Proof Assertion (i) is evident since $\Gl>0$. Assertion (ii), as well as the existence part of 
assertion (iii), can be found in \cite {HJV} or \cite {Ve5}. Furthermore any $\gw\in \CS^*_{q\,+}$ 
is positive in $S^{N-1}_{+}$ and verifies $\gw_{\gf}<0$ by Hopf boundary lemma as the 
outward normal derivative on $\prt S^{N-1}_{+}$ is $\prt\,/\prt\gf$.  We can construct a 
minimal element in  
$\CS^*_{q\,+}$ in the following way: If we denote by $u^H_{k}$ the unique solution of 
(\ref {equ+abs}) in $H$ which satisfies (\ref{sing-k}) and set $T_{r}(u^H_{k})(x)
=r^{\gb_{q}}u^H_{k}(rx)$ for $r>0$, then $T_{r}(u^H_{k})$ is a solution of (\ref {equ+abs}) 
in $H$ which satisfies
$$T_{r}(u^H_{k})=r^{(2N-1-q)/(q+1-N)}k\myfrac {x_{N}}{\abs x^2}(1+\circ (1))\quad \mbox 
{as }x\to 0.
$$
Thus $T_{r}(u^H_{k})=u^H_{r^{(2N-1-q)/(q+1-N)}k}$. Furthermore, if $\gw\in\CS^*_{q\,+}$, the 
maximum principle at $0$ and at infinity (replacing $u_{\gw}$ by $u_{\gw}+\ge$ and letting $\ge\to 0$) leads to 
$$u_{\gw}(r,\gs):=r^{-\gb_{q}}\gw(\gs)>u^H_{k}(r,\gs)\forevery (r,\gs)\in (0,\infty)\ti 
S^{N-1}_{+},\,\forall k>0.
$$
Letting $k\to\infty$ implies $u_{\gw}(r,\gs)\geq u^H_{\infty}(r,\gs)$ and 
$T_{r}(u^H_{\infty})=u^H_{\infty}$ given that $2N-1-q>0$. Then the function $u^H_{\infty}$ is invariant 
with respect to the transformation $T_{r}$. It is therefore self-similar, and consequently under 
the form
$u^H_{\infty}(r,\gs)=r^{-\gb_{q}}\underline\gw(\gs)$. As a result of the previous inequality 
$\underline\gw$ is the minimal element of $\CS^*_{q\,+}$. Next we denote 
$\gd^*=\max\{\gd\geq 0:\gd\gw\leq\underline\gw\}$ and $u_{\gw,\gd^*}=\gd^* u_{\gw}$. 
Notice that $\gd^*\in (0,1]$ as $\underline\gw>0$ in $S^{N-1}_{+}$ and satisfies Hopf 
boundary lemma on $\prt S^{N-1}_{+}$
Clearly $u_{\gw,\gd^*}$ is a subsolution for (\ref{equ+abs}) and it is dominated by $u^H_{\infty}$ 
in $H$. Furthermore $\gd^*\gw\leq\underline\gw$ in $\overline {S^{N-1}_{+}}$, 
$\gd^*\gw_{\gf}\leq\underline\gw_{\gf}$ on $\prt S^{N-1}_{+}$, and \smallskip

\nind (i) either there exists $\gs_{0}\in S^{N-1}_{+}$ such that $\gd^*\gw(\gs_{0})
=\underline\gw(\gs_{0})$,  \smallskip

\nind  (ii) or $\gd^*\gw<\underline\gw$ in $S^{N-1}_{+}$ and there exists $\gs'_{0}\in S^{N-2}$ 
such that
$\gd^*\gw_{\gf}(\gs'_{0},\gp/2)=\underline\gw_{\gf}(\gs'_{0},\gp/2)$.\smallskip

\nind In case (i), and as $Du^H_{\infty}$ never vanishes in $H$, it follows from 
\cite[Lemma 1.3]{FV}
(a variant of the strong comparison principle) that $u_{\gw,\gd^*}=\underline u$. This implies that 
$u_{\gw,\gd^*}$ is a solution, $\gd^*=1$  and, consequently $\gw=\underline\gw$.\smallskip

\nind In case (ii) we follow the linearization procedure already introduced in \cite {FV}. 
By the mean value theorem 
$$\abs {Du^H_{\infty}}^{N-2}u_{\infty\,x_{i}}-
\abs {Du_{\gw,\gd^*}}^{N-2}u_{\gw,\gd\,x_{i}}=\mysum{j}{}\ga_{ij}
(u^H_{\infty}-u_{\gw,\gd^*})_{x_{j}}
$$
where
$$\BA {l}
\ga_{ij}=\abs{t_{i}Du^H_{\infty}+(1-t_{i})Du_{\gw,\gd^*}}^{N-4}\left(\gd_{ij}\abs{t_{i}
Du^H_{\infty}+(1-t_{i})Du_{\gw,\gd^*}}^{2}\right.\\[2mm]
\phantom {--------}
\left.+(N-2)
\left(t_{i}u^H_{\infty\,x_{i}}+(1-t_{i})u_{\gw,\gd^*\,x_{i}}\right)
\left(t_{i}u^H_{\infty\,x_{j}}+(1-t_{i})u_{\gw,\gd^*\,x_{j}}\right)\right),
\EA$$
with $0\leq t_{i}\leq 1$. Next $w=u^H_{\infty}-u_{\gw,\gd^*}$ is positive in $H$ and satisfies
$$-\mysum{ij}{}\left(\ga_{ij}w_{x_{j}}\right)_{x_{i}}+cw\geq 0
$$
where $c=((u^H_{\infty})^q-u^q_{\gw,\gd^*})/(u^H_{\infty}-u_{\gw,\gd^*})>0$.
Notice that $(\ga_{ij}(x))$ is the Hessian of a strictly convex function  therefore 
it is nonnegative and that $(\ga_{ij})(r,\gs'_{0},\gp/2)$ is positive-definite. Therefore 
it is positive-definite in a neighborhood of 
$(r,\gs'_{0},\gp/2)$ (independent of $r$, actually). 
Inasmuch $(u^H_{\infty}-u_{\gw,\gd^*})_{x_{N}}=0$ at 
$(r,\gs'_{0},\gp/2)$, we derive a contradiction with Hopf lemma. Therefore case (ii) cannot 
occur and 
$\gw=\underline\gw$.\qeda
\\

 %%%%%%%%%%%%%%%%%%%%%%%%%%%%%%%%%%%%%%%%%%%%%%%
%%%%IMPORTANT REMARK %%%%%%%%%%%%%%%%%%%%%%%%%%
%%%%%%%%%%%%%%%%%%%%%%%%%%%%%%%%%%%%%%%%%%%%%%%
\nind \Remark If we look for separable solutions of
\begin {equation}\label {equ+abs-p}
-div\left(\abs {Du}^{p-2}Du\right)+\abs u^{q-1}u=0,
\end {equation}
in $\BBR^N$, where $q>p-1>0$, $p$ not necessarily equal to $N$ or to $2$, under the form 
$u(r,\gs)=r^{\gb}\gw(\gs)$, then $\gb=\gb_{p,q}=-p/(q+1-p)$ and $\gw$ is a solution of 
\begin {equation}\label {sep2}
-div_{\gs}\left(\left(\gb^2_{p,q}\gw^2+
\abs {\nabla_{\gs}\gw}^{2}\right)^{(p-2)/2}\nabla_{\gs}\gw\right)
-\Gl(p,q)\left(\gb^2_{p,q}\gw^2+\abs {\nabla_{\gs}\gw}^{2}\right)^{(p-2)/2}\gw
+\abs w^{q-1}\gw=0
\end {equation}
on $S^{N-1}$ where $\Gl(p,q)=\gb^{p-1}_{p,q}(q\gb_{p,q}-p)$. If we look for 
separable solutions in $H$ which vanishes on $\prt H\setminus \{0\}$ the solution 
$\gw$ of (\ref {sep2}) is subject to 
the boundary condition $\gw=0$ on $\prt S^{N-1}_{+}=S^{N-2}$. A fairly exhaustive theory of 
existence is developped in \cite {Ve5}, \cite {HJV}. The existence of non-trivial solution 
of (\ref {sep2}) is insured as soon $\Gl(p,q)>0$, or equivalently $q<N(p-1)/(N-p)$ if $p<N$, 
and no condition if $p\geq N$. If $q\geq N(p-1)/(N-p)$ no solution exists, up to the trivial one. 
This is linked to the removability result proved by V\`azquez and V\'eron \cite {VV}. The existence 
of non trivial solutions of the same equation in $S^{N-1}_{+}$ vanishing on $\prt S^{N-1}_{+}$ 
is much more complicated. However it is proved in \cite {Ve5}, \cite {HJV} that there exists a 
critical exponent $q_{c}>p-1$ such that, if $q\geq q_{c}$ no non-trivial solution exists while 
if $p-1<q<q_{c}$ there exist a unique positive solution in $S^{N-1}_{+}$ vanishing on 
$\prt S^{N-1}_{+}$. The uniqueness proof in the previous proposition is valid. \\

The next result characterizes the solution of (\ref {equ+abs}) with a strong singularity 
on the boundary. In order to express the result, we assume that the outward normal 
unit vector to $\prt\Gw$ at $a$ is $-{\bf e}_{N}$. 
 %%%%%%%%%%%%%%%%%%%%%%%%%%%%%%%%%%%%%%%%%%%%%%%
%%%%BIG THEOREM II %%%%%%%%%%%%%%%%%%%%%%%%%%
%%%%%%%%%%%%%%%%%%%%%%%%%%%%%%%%%%%%%%%%%%%%%%%
\bth {gendom2} Let $\Gw$ be a bounded domain with a $C^2$ boundary and $a\in\prt\Gw$. 
Assume $0<p-1<q<2N-1$. Then  there exists a unique function 
$u\in C^1(\overline \Gw\setminus \{a\})$ which vanishes on $\prt\Gw\setminus\{a\}$, 
is solution of (\ref {equ+abs}) in $\Gw$ and satisfies
\begin {equation}\label {sing-k10}
 \lim_{x\to a}\myfrac {\abs {x-a}^{2}u(x)}{\gr(x)}=\infty.
\end {equation}
Furthermore 
\begin {equation}\label {sing-k11}
 \lim_{\scz {\BA {c}x\to a\\
 (x-a)/\abs {x-a}\to\gs\EA}}\!\!\!\!\!\!\abs {x-a}^{\gb_{q}}u(x)=\gw(\gs),
\end {equation}
locally uniformly on $S^{N-1}_+$. Finally $u=u_{\infty,a}=\lim_{k\to\infty}u_{k,a}$.
\es
\Proof We already know that $u_{\infty,a}$ satisfies (\ref{sing-k10}). 
By translation we fix the origin $0$ of coordinates at the point $a$ 
and we assume that $-{\bf e}_{N}$ is the outward unit vector to 
$\prt\Gw$ at $0$.
If $G$ is any $C^2$ domain in $\BBR^N$ to the boundary of which $0$ belongs, 
we denote by $u^{G}_{k}$ the solution of (\ref {equ+abs}) in $G$, 
which vanishes on $\prt G\setminus\{0\}$ and verifies 
\begin {equation}\label {sing-k12}
u^{G}_{k}=k\myfrac {\gr_{_G}(x)}{\abs x^{2}}(1+\circ (1))\quad\mbox { as }x\to 0,
\end {equation}
where $\gr_{_G}(x)=\dist(x,G)$. When there is no ambiguity, $u^\Gw_{k}=u_{k}$. By the maximum principle 
$G\subset G'$ implies $u^{G}_{k}\leq u^{G'}_{k}$ 
in $G$. By dilation we can assume that there exist two balls of 
radius $1$, $B\subset \Gw$ and $B'\subset \overline {\Gw}^c$ with 
respective center $b={\bf e}_{N}$ and $b'=-b$ with the property that 
$0=\prt B\cap \prt B'$. It follows from the maximum principle, the 
fact that $u_{k}^{B}(x)=u_{k}^{B'}(\CS(x))$ where $\CS$ is the 
symmetry with respect to the hyperplane $\prt H$ and \rprop {ball}, \rprop 
{compball}
\begin {equation}\label {comp1}\BA {l}
(i)\quad u_{k}^{B}(x)\leq u_{k}(x)\leq u_{k}^{B'\,\!^c}(x)\leq 
u_{k}^{B'}\left(b'+\myfrac{x-b'}{\abs {x-b'}^{2}}\right)=
u_{k}^{B}\left(\CS\left(b'+\myfrac{x-b'}{\abs {x-b'}^{2}}\right)\right)\\
\phantom {u_{k}^{B}(x)\leq u_{k}(x)\leq u_{k}^{B'\,\!^c}(x)\leq 
u_{k}^{B'}\left(b'+\myfrac{x-b'}{\abs {x-b'}^{2}}\right)=---------}
\forevery x\in 
B\\[2mm]
(ii)\quad u_{k}(x)\leq u_{k}^{B'\,\!^c}(x)\leq 
u_{k}^{B}\left(\CS\left(b'+\myfrac{x-b'}{\abs {x-b'}^{2}}\right)\right)\forevery x\in \Gw,
\EA
\end {equation}
and similarly
\begin {equation}\label {comp2}\BA {l}
(i)\quad u_{k}^{B}(x)\leq u^{H}_{k}(x)\leq u_{k}^{B'\,\!^c}(x)\leq 
u_{k}^{B}\left(\CS\left(b'+\myfrac{x-b'}{\abs {x-b'}^{2}}\right)\right)\forevery x\in 
B\\[4mm]
(ii)\quad u^{H}_{k}(x)\leq u_{k}^{B'\,\!^c}(x)\leq 
u_{k}^{B}\left(\CS\left(b'+\myfrac{x-b'}{\abs {x-b'}^{2}}\right)\right)\forevery x\in 
H,
\EA
\end {equation}
Letting $k\to\infty$, we obtain
\begin {equation}\label {comp3}\BA {l}
(i)\quad u_{\infty}^{B}(x)\leq u_{\infty}(x)\leq 
u_{\infty}^{B'\,\!^c}(x)\leq 
u_{\infty}^{B}\left(\CS\left(b'+\myfrac{x-b'}{\abs {x-b'}^{2}}\right)\right)\forevery x\in 
B\\[4mm]
(ii)\quad u_{\infty}(x)\leq u_{\infty}^{B'\,\!^c}(x)\leq 
u_{\infty}^{B}\left(\CS\left(b'+\myfrac{x-b'}{\abs {x-b'}^{2}}\right)\right)\forevery x\in \Gw,
\EA\end {equation}
as well as
\begin {equation}\label {comp4}\BA {l}
(i)\quad u_{\infty}^{B}(x)\leq 
\abs x^{-\gb_{q}}\gw\left(\myfrac{x}{\abs x}\right)
\leq u_{\infty}^{B'\,\!^c}(x)\leq 
u_{\infty}^{B}\left(\CS\left(b'+\myfrac{x-b'}{\abs {x-b'}^{2}}\right)\right)\forevery x\in 
B\\[4mm]
(ii)\quad \abs x^{-\gb_{q}}\gw\left(\myfrac{x}{\abs x}\right)\leq 
u_{\infty}^{B'\,\!^c}(x)
\leq 
u_{\infty}^{B}\left(\CS\left(b'+\myfrac{x-b'}{\abs {x-b'}^{2}}\right)\right)\forevery x\in H.
\EA
\end {equation}
From (\ref {comp4})-(i) and the  fact that $b'=-b$, we also derive
\begin {equation}\label{comp5}\BA {l}
\abs x^{-\gb_{q}}\gw(x/\abs x)\leq 
u_{\infty}^{B'\,\!^c}(x)
\leq 
\abs {\CS\left(\myfrac{x+b}{\abs {x+b}^{2}}-b\right)}^{-\gb_{q}}\!\!\!\!\!
\gw\left(\myfrac {\CS\left(x+b-\abs {x+b}^{2}b\right)}{\abs{\CS\left(x+b-\abs {x+b}^{2}b\right)}}\right).
\EA\end {equation}
But 
$$\abs {\CS\left(\myfrac{x+b}{\abs {x+b}^{2}}-b\right)}=\myfrac {\abs x}{\abs {x+b}}=
\abs x(1+\circ(1))\quad\mbox {as }x\to 
0$$ 
(remember that $\abs b=1$). If $x=(x_{1},\ldots,x_{N})$, $\abs 
{x+b}^{2}=\abs x^{2}+1+2x_{N}$ and
$$\CS\left(x+b-\abs {x+b}^{2}b\right)=(x_{1},\ldots,x_{N}+\abs x^{2}).$$
%%%%%%%%%%%%%%%%%%%%%%%%%%%
Thus (\ref{comp5}) becomes 
\begin {equation}\label{comp6}\BA {l}
\abs x^{-\gb_{q}}\gw(x/\abs x)\leq 
u_{\infty}^{B'\,\!^c}(x)
\leq \abs x^{-\gb_{q}}\abs {x+b}^{\gb_{q}}\gw\left(\myfrac {x+\abs x^{2}{\bf e}_{N}}{\abs 
x\sqrt {1+\abs x^{2}+2x_{N}}}\right).
\EA
\end {equation}
If we assume $\abs x^2=\circ(x_{N})$ then $(x+\abs x^{2}{\bf e}_{N})/(\abs 
x\sqrt {1+\abs x^{2}+2x_{N}})=x(1+\circ (1))/\abs x$ as $x\to 0$, and
\begin {equation}\label{comp7}\BA {l}
u_{\infty}^{B'\,\!^c}(x)=\abs x^{-\gb_{q}}\gw(x/\abs x)(1+\circ (1)).
\EA
\end {equation}
If we define $\CT$ by 
$$\CT(x)=\CS\left(\myfrac{x+b}{\abs {x+b}^{2}}-b\right),
$$
then (\ref {comp4})-(i) reads also as 
\begin {equation}\label {comp8}\BA {l}
\abs {\CT^{-1}(x)}^{-\gb_{q}}\gw\left(\myfrac{\CT^{-1}(x)}{\abs {\CT^{-1}(x)}}\right)\leq  u_{\infty}^{B}(x)\leq 
\abs x^{-\gb_{q}}\gw\left(\myfrac{x}{\abs x}\right).
\EA
\end {equation}
Furthermore
$$\CT^{-1}(x)=\left(\myfrac {x_{1}}{\abs {x-b}^2},...,\myfrac {x_{N-1}}{\abs {x-b}^2},\myfrac {1-x_{N}}{\abs {x-b}^2}-1\right)=\myfrac {x-\abs x^2{\bf e}_{N}}{\abs {x-b}^2}.
$$
Then
$$\abs{\CT^{-1}(x)}=\abs{b+\myfrac {x-b}{\abs {x-b}^2}}=\myfrac {\abs x}{\abs {x-b}},
$$
and 
$$\abs {\CT^{-1}(x)}^{-\gb_{q}}\gw\left(\myfrac{\CT^{-1}(x)}{\abs {\CT^{-1}(x)}}\right)
=\abs x^{-\gb_{q}}\abs {x-b}^{\gb_{q}}\gw\left(\myfrac {x-\abs x^2{\bf e}_{N}}{\abs x\abs {x-b}}\right).
$$
If we assume again $\abs x^2=\circ(x_{N})$ then $(x-\abs x^{2}{\bf e}_{N})/(\abs 
x\sqrt {1+\abs x^{2}-2x_{N}})=x(1+\circ (1))/\abs x$ as $x\to 0$, and
\begin {equation}\label{comp9}\BA {l}
u_{\infty}^{B}(x)=\abs x^{-\gb_{q}}\gw(x/\abs x)(1+\circ (1)).
\EA
\end {equation}
Combining (\ref{comp3})-(i), (\ref{comp6}) and (\ref{comp8}) we obtain that
\begin {equation}\label{comp10}\BA {l}
u_{\infty}(x)=\abs x^{-\gb_{q}}\gw(x/\abs x)(1+\circ (1))\quad \mbox { as }x\to 0
\EA
\end {equation}
uniformly on any subset of $\Gw$ such that $\abs x^2=\circ(x_{N})$ near $0$. In order to obtain the precise behaviour (\ref {sing-k11}), we proceed and in the proof of \rth {gendom}. We extend $u$ by reflection through $\prt\Gw$ near $0$ and denote by $\tilde u$ the extended function defined in $T^\gb(\Gw)$. For $r\in (0,1]$ we define
$$w_{r}:=T_{r} (\tilde u)(x)=r^{\gb_{q}}\tilde u(r x).
$$
Then $w_{r}$ satisfies
  \begin {equation}\label {ext-k3}
-\mysum {j}{}\myfrac {\prt }{\prt x_{j}}\tilde A^r_{j}(x,Dw_{r})+\tilde b(rx)\abs {w_{r}}^{q-1}w_{r}=0
\end {equation}
 in $T^{\gb r^{-1}}(\Gw^{r})$. Since 
$w_{r}$ is uniformly bounded on $\Gg_{1/2,2}$ (by \rprop {KOV} applied to $u$ and $-u$) and the definition of the refected function), $Dw_{r}(u)$ is bounded in $C^{\ga}(\Gg_{2/3,3/2})$. By Ascoli's theorem $w_{r}$ converges in the $C^1(\Gg_{2/3,3/2})$-topology to $x\mapsto \abs x^{-\gb_{q}}\tilde \gw(x/\abs x)$, where $\tilde \gw$ is defined from $\gw$ by reflection through the equator $\prt S^{N-1}_{+}$. In ordre to get rid of the boundary, we use again the $C^2$ diffeomorphism $\Psi$ which sends $B_{\ge_0}$ onto itself and verifies $\Psi(B_{\ge_0}\cap\prt\Gw)=B_{\ge_0}\cap\prt H$. We  set $x=\Psi^{-1}(y)$ and $\tilde u(x)=u^*(y)$. Then
\begin {equation}\label {ext-k4}
\lim_{\abs y\to 0} (\abs y^{\gb_{q}} u^*(\abs y,\gs),\abs y^{\gb_{q}+1}\nabla u^*(\abs y,\gs))=
(\gw(\gf),-\gw_{\gf} \,\bf{e}_N ),
\end {equation}
uniformly on $S^{N-1}$, thus
\begin {equation}\label {asym''3}
 u^*(\abs y,\gs)=\abs y^{\gb_{q}}\gw(\gf) (1+\circ (1))\quad \mbox {as }\,\abs y\to 0
\end {equation}
uniformly with respect to $\gs\in S^{N-1}_+$, because $u^*$ vanishes on $B_{\ge_0}\cap\prt H\setminus\{0\}$. Actually, a stronger result than (\ref{sing-k11}) follows, namely
  \begin {equation}\label {sing-k13}
u(x)=\abs x^{-\gb_{q}}\gw(x/\abs x)(1+\circ 1))\quad\mbox {as }x\to 0.
\end {equation}
{\it Mutatis mutandis}, this estimate implies uniqueness of a solution with a strong singularity
 as in \rth {gendom}.\qeda
%%%%%%%%%%%%%%%%%%%%%%%%%%%%%%%%%%%%%%%%%%%%%%%
%%%%%%%%%SECTION REMOVABILITY%%%%%%%%%%%%%%%%%%%%%%%%%%%%%
%%%%%%%%%%%%%%%%%%%%%%%%%%%%%%%%%%%%%%%%%%%%%%%%%%%%
\mysection {The removability result}
In this section $\Gw$ is a $C^2$ domain of $\BBR^N$ and $a\in\prt\Gw$. The next result extends the removability result of Gmira-V\'eron \cite {GV} dealing with semilinear equations.
%%%%%%%%%%%%%%%%%%%%%%%%%%%%%%%%%%%%%%%%%%%%%%%
%%%%%%%%%THEOREM REMOVABILITY%%%%%%%%%%%%%%%%%%%%%%%%%%%%%
%%%%%%%%%%%%%%%%%%%%%%%%%%%%%%%%%%%%%%%%%%%%%%%%%%%%
\bth {remov}Let $g$ be a continuous function defined on $\BBR$ which satisfies
\begin{equation}\label {cond-g}
\liminf_{r\to\infty}g(r)/r^{q_{c}}>0\quad\mbox {and }\;\;\limsup_{r\to -\infty}g(r)/|r|^{q_{c}}<0,
\end{equation}
where $q_{c}:= 2N-1$ and let $u\in C^1(\overline\Gw\setminus\{a\})$ be a solution of 
\begin{equation}\label {remov'1}
-div\left(\abs {Du}^{N-2}Du\right)+g(u)=0\quad\mbox {in }\Gw
\end{equation}
which coincides with some $\gf\in C^1({\prt\Gw})$ on $\prt\Gw\setminus\{a\}$. Then $u$ extends to $\overline\Gw$ as a continuous function.
\es
\Proof Without any loss of generality, we can assume that $\Gw$ is bounded, $a=0$ and $-{\bf e}_{N}$ is the outward normal vector to $\prt\Gw$ at $0$. We denote by $V_{0}^\Gw$ the solution of (\ref{equ+hom}) in $\Gw$ which vanishes on $\prt\Gw\setminus\{0\}$ and satisfies
$$V_{0}^\Gw(x)=\myfrac {\gr(x)}{\abs {x}^2}(1+\circ (1))\quad\mbox {as }x\to 0.
$$
Let $M$ be the supremum of $\abs\gf$ on $\prt\Gw$ and $\tilde M=\max\{M,\left(B/A\right)^{1/q}\}$. By assumption there exists $A>0$ and $B\geq 0$, depending only on $g$, such that
\begin{equation}\label {remov2}
-div\left(\abs {Du}^{N-2}Du\right)+A u^{q_{c}}\leq B\quad\mbox {in }\{x\in\Gw:u(x)> 0\}.
\end{equation}
If $ v=u-\tilde M$, then $v\leq 0$ on $\prt\Gw\setminus\{0\}$ and 
\begin{equation}\label {remov3}
-div\left(\abs {Dv}^{N-2}Dv\right)+A v^{q_{c}}\leq 0\quad\mbox {in }\{x\in\Gw: v(x)> 0\},
\end{equation}
Using the same functions $\eta_{\ge}$ as in the proof of \rprop{KOV} we deduce that 
$\eta_{\ge}(v)$ satisfies the same inequality as $v$, but on whole $\Gw$. 
By \rprop {UPREG} with $q=q_{c}$ and the expression of $V_{0}^\Gw$ it follows that 
\begin{equation}\label {remov4}
v(x)\leq cV^\Gw(x)\forevery x\in\Gw,
\end{equation}
where the constant $c$ depends on $A$ and $N$. Furthermore, there exists a function $u^*\in C^1(\overline\Gw\setminus \{0\})$ such that $0\leq v_{+}\leq u^*(x) \leq cV^\Gw_{0}$ in $\Gw$, and
\begin{equation}\label {remov5}
-div\left(\abs {Du^*}^{N-2}Du^*\right)+A u^{*\,q_{c}}=0\quad\mbox {in }\Gw.
\end{equation}
As in the proof of \rth {gendom2} we extend $u^*$ through the boundary into $\tilde u$ and scale it by setting  $T_{r}(\tilde u):=w_{r}(x)=r\tilde u(rx)$ for $r>0$. Inasmuch all the previous a priori estimates apply (compactness), it follows that 
there exists a subsequence $\{r_{n}\}$ converging to $0$ and a function $w\in C^1(\BBR^N\setminus\{0\})$ such that $w_{r_{n}}\to w$ in the $C^1_{loc}$-topology of $\BBR^N\setminus\{0\}$, $w$ is a solution of 
\begin{equation}\label {remov6}\left\{\BA {l}
-div\left(\abs {Dw}^{N-2}Dw\right)+A w^{q_{c}}=0\quad\mbox {in }\BBR^N\setminus\{0\}\\
\phantom{div\left(\abs {Dw}^{N-2}Dw\right)+A w^{q_{c}}}
w\geq 0\;\mbox { in }H=\{x\in\BBR^N:x_{N}>0\}\\
\phantom{div\left(\abs {Dw}^{N-2}Dw\right)+A w^{q_{c}}}
w=0\;\mbox { on }\prt H\setminus\{0\}.
\EA\right.
\end{equation}
At end, (\ref{remov4}) transforms into 
\begin{equation}\label {remov7}
0\leq w(x)\leq c\myfrac{x_{N}}{\abs x^2}\forevery x\in H.
\end{equation}
For $\ge>0$ we denote by $W_{\ge}$ the solution of 
\begin{equation}\label {remov8}\left\{\BA {l}
-div\left(\abs {DW_{\ge}}^{N-2}DW_{\ge}\right)+A W_{\ge}^{q_{c}}=0\quad\mbox {in }H\setminus B_{\ge}(0)\\
\phantom{div\left(\abs {DW_{\ge}}^{N-2}DW_{\ge}\right)+A W_{\ge}}
W_{\ge}=c\ge^{-2}x_{N}\;\mbox { on }H\cap \prt B_{\ge}(0)\\
\phantom{div\left(\abs {DW_{\ge}}^{N-2}DW_{\ge}\right)+A W_{\ge}}
W_{\ge}=0\;\mbox { on }\prt H\setminus B_{\ge}(0).
\EA\right.
\end{equation}
By the maximum principle $0\leq w(x)\leq W_{\ge}(x)\leq cx_{N}\abs x^{-2}$ for any $\ge>0$, and by uniqueness, 
$T_{r}(W_{\ge})(x)=rW_{\ge}(rx)=W_{\ge/r}(x)$. Furthermore $\ge\mapsto W_{\ge}$ is increasing. Letting $\ge\to 0$ we conclude that $W_{\ge}$ decreases to some $W_{0}$, which is a solution of 
\begin{equation}\label {remov'8}\left\{\BA {l}
\!\!-div\left(\abs {DW_{0}}^{N-2}DW_{0}\right)+A W_{0}^{q_{c}}=0\quad\mbox {in }H\\
\phantom{div\left(\abs {DW_{\ge}}^{N.2}DW\right)+A W_{\ge}^{q_{c}}}
W_{0}\geq 0\;\mbox { in }H\\
\phantom{div\left(\abs {DW_{\ge}}^{N.2}DW\right)+A W_{\ge}^{q_{c}}}
W_{0}=0\;\mbox { on }\prt H\setminus \{0\},
\EA\right.
\end{equation}
by the standard regularity results, and satisfies $0\leq w\leq W_{0}$. Finally, $W_{0}$ inherits the following scaling invariance property $T_{r}(W_{0})(x)=W_{0}(x)$ for any $r>0$. Therefore $W_{0}$ is a separable solution which endows the following form
$$W_{0}(x)=W_{0}(r,\gs)=r^{-1}\gw(\gs),
$$
where $\gw$ is nonnegative on $S^{N-1}_{+}$ and satisfies
\begin {equation}\label {remov9}\left\{\BA {l}
\!\!\!-div_{\gs}\left(\left(\gw^2+\abs {\nabla_{\gs}\gw}^{2}\right)^{(N-2)/2}\!\!\!
\nabla_{\gs}\gw\right)-(N-1)\left(\gw^2+\abs {\nabla_{\gs}\gw}^{2}
\right)^{(N-2)/2}\!\!\!\!\!\gw+ A\gw^{q_{c}}=0\\
\phantom{-----------------------------------}
\mbox {in }S^{N-1}_{+}\\[2mm]
\phantom{------------------------------}
\gw=0\quad\,\mbox { on }\prt S^{N-1}_{+}.\\
\EA\right.\end {equation}
By \rprop {psirad}, $\gw=0$. Thus $W_{0}=0\Longrightarrow w=0$, which implies  $w_{r}(x)\to 0$ as $r\to 0$
and equivalently $r\tilde u(rx)\to 0$ in the $C^1_{loc}$-topology of $\BBR^N\setminus\{0\}$. Consequently
$D\tilde u(x)=\circ (\abs x^{-2})$ as $x\to 0$ and finally $u^*(x)=\circ (V_{0}^\Gw(x))$ as $x\to 0$. The maximum principle and the positivity of $u^*$ yields to $u^*\equiv 0$ and finally $u\leq \tilde M$ in $\Gw$. In the same way $u\geq -\tilde M$. Because the modulus of continuity of $u$ is uniformly bounded near $0$, by the classical regularity theory of degenerate elliptic equations (see \cite {Lie} for example), $u$ extends as a continuous function in whole
$\overline\Gw$.\qeda
\section {The classification theorem}

The next result extends some of Gmira-V\'eron's classification theorem \cite [Sect. 4, 5 ]{GV} obtained in the study of problem (\ref{lin2}). In the above mentioned article, the main idea was to reduce the equation to a infinite dimensional quasi-autonomous evolution system in $\BBR_{+}\ti S^{N-1}_{+}$ and to use Lyapounov-energy function. Such an approach cannot be adapted in the quasilinear case. Our method is based upon scaling and uniqueness arguments.
%%%%%%%%%%%%%%%%%%%%%%%%%%%%%%%%%%%%%%%%%%%%%%%%%%%%%
%%%%%%%%%%%%% THEOREM-CLASS 1%%%%%%%%%%%%%%%%%%%%%%%%%%%%%%%%%%%%%%%%%%%%%%%%%%%%%%%%%%%%%%%%%%%%%%%%%%%%%%%%%%%%
\bth {class}Assume $N-1<q<2N-1$, $\Gw$ is a bounded domain with a $C^2$ boundary, $a\in\prt\Gw$ and $-{\bf e}_{N}$ is the outward normal unit vector to $\prt\Gw$ at $a$. Let $u\in C^1(\overline\Gw\setminus\{a\})$ be a positive function satisfying (\ref{equ+abs}) in $\Gw$ and vanishing on $\prt\Gw\setminus\{a\}$. Then the following alternative holds.\smallskip

\nind (i) Either there exists $k>0$ such that 
\begin{equation}\label{sing1}
u(x)=k\myfrac{\gr(x)}{\abs{x-a}^2}(1+\circ (1))\quad\mbox {as }x\to a.
\end {equation}
Furthermore $u=u_{k,a}$, the unique solution of (\ref {equ+abs}) defined in \rth {gendom}.\smallskip

\nind (ii) Or 
\begin{equation}\label{sing2}
u(x)=\abs {x-a}^{-\gb_{q}}\gw(\gs)(1+\circ (1))\quad\mbox {as }\,x\to a.
\end {equation}
where $\gw$ is the unique positive solution of (\ref{sep1}) on $S^{N-1}_{+}$ which vanishes on $\prt S^{N-1}_{+}$, in which case $u=u_{\infty,a}$.
\es
%%%%%%%%%%%%%%%%%%%%%%%%%%%%%%%%%%%%%%%%%%%%%%%%%%%%%
%%%%%%%%%%%%%PROOF OF THEOREM-CLASS 1%%%%%%%%%%%%%%%%%%%%%%%%%%%%%%%%%%%%%%%%%%%%%%%%%%%%%%%%%%%%%%%%%%%%%%%%%%%%%%%%%%%%
\Proof We assume that $a=0$ with $\gn_{0}=-\bf {e}_{N}$ and define
\begin{equation}\label{S1}
k=\limsup_{x\to 0}\myfrac{u(x)}{V_{0}^{\Gw}(x)}
=\limsup_{r\to 0}\sup_{\abs x=r}\myfrac{u(x)}{V_{0}^{\Gw}(x)}.
\end{equation}
Suppose $k=0$. It follows from the maximum principle that for any $\ge>0$ there exists a sequence $r_{n}\to 0$ such that $0\leq u(x)\leq \ge V_{0}^{\Gw}(x)$ in $\Gw\setminus\{B_{r_{n}}(0)\}$. This fact implies the nullity of $u$. Therefore we assume that $k\neq 0$. Assume first that $k$ is finite. Then, for any $\ge>0$, there exists a sequence of points $x_{n}$ converging to $0$ such that 
\begin{equation}\label{S2}
\lim_{n\to \infty} \myfrac{u(x_{n})}{V_{0}^{\Gw}(x_{n})}=k
\end{equation}
and
\begin{equation}\label{S3}
\sup_{\abs x\leq r_{n}} \myfrac{u(x)}{V_{0}^{\Gw}(x)}\leq k+\ge.
\end{equation}
Since $u_{k}$ satisfies (\ref{sing-k3}) with $a=0$, the two previous relations can be replaced by
\begin{equation}\label{S2'}\BA {l}
(i)\phantom{-----}\,\lim_{n\to \infty} \myfrac{u(x_{n})}{u_{k}(x_{n})}=1\phantom{------------------}\\
(ii)\phantom{-----}\sup_{\abs x\leq r_{n}} \myfrac{u(x)}{u_{k}(x)}\leq 1+\ge.\phantom{------------------}
\EA\end{equation}
We denote $r_{n}=\abs {x_{n}}$, $\xi_{n}=x_{n}/r_{n}$ and define $u_{n}=r_{n}u(r_{n}x)$ and
$u_{k\,n}=r_{n}u_{k}(r_{n}x)$. By the previous arguments combining a priori estimate and regularity theory, there exist a subsequence $\{r_{n_{j}}\}$ and two nonnegative functions $v$ and $v'$,  $N$-harmonic in $H$ and vanishing on $\prt H\setminus\{0\}$, such that 
$\left(u_{n_{j}},u_{k\,n_{j}}\right)$ converges to $(v,v')$ in the $C^1_{loc}$-topology of $H=\BBR^N_{+}$. Clearly equality (\ref{sing-k3}) implies that $r_{n_{j}}u_{k}(r_{n_{j}}x)$ converges to $kV_{0}^{H}(x)$ (which is defined by $kV_{0}^{H}(x):= kx_{N}/\abs x^{-2}$) in the same topology. Since $v'$ is uniquely determined by its blow-up at $0$, this  implies $v'=kV_{0}^{H}$ in $H$.
Furthermore there exists $\xi\in \overline{S^{N-1}_{+}}$  such that $\xi_{n_{k}}\to\xi$. If $\xi\in S^{N-1}_{+}$, $v(\xi)=v'(\xi)$, while, if $\xi\in \prt S^{N-1}_{+}$, 
$\prt v/\prt \gn(\xi)=\prt _{x_{N}}v(\xi)=\prt v'/\prt \gn(\xi)$. In both situation, the tangency conditions of the graphs of $v$ and $v'$ and the strong maximum principle implies that 
$v=v'=kV_{0}^{H}$. By estimate (\ref {S2'})-(i) and the convergence properties, it follows
$$\lim_{n\to \infty} \myfrac{u(r_{n}\xi)}{u_{k}(r_{n}\xi)}=1,\quad \mbox {uniformly on }\abs \xi=1.
$$
Consequently, for any $\gd>0$, there holds,
$$(1-\gd)u_{k}(x)\leq u(x)\leq (1+\gd)u_{k}(x)\forevery x\in\Gw\setminus B_{r_{n}},
$$
 for $n$ large enough, which leads to $u_{k}= u$. At end we consider the case $k=\infty$. Writting (\ref{equ+abs}) under the form
\begin{equation}\label{sing3}
-div\left(\abs {Du}^{N-2}Du\right)+d(x)u^{N-1}=0
\end {equation}
where $d(x)=\abs {u}^{q+1-N}(x)\leq C\abs x^{-N}$, by (\ref{up-est}), We use the boundary Harnack principle. By \cite [Th 2.2]{BBV} there exists a constant $c=c(N,q,\Gw)>0$ such that 
\begin{equation}\label{sing4}
\myfrac{1}{c}\myfrac{u(y)}{\gr(y)}\leq \myfrac{u(x)}{\gr(x)}\leq c\myfrac{u(y)}{\gr(y)}
\end {equation}
for any $x$ and $y$ in $\Gw$ such that $\abs x=\abs y$ be small enough. Since there exists a sequence $x_{n}\to 0$ such that $\lim_{n\to\infty}u(x_{n})/V^{\Gw}(x_{n)}\to\infty$, this implies that
\begin{equation}\label{sing5}
\lim_{n\to\infty}\left\{\inf\myfrac{u(x)}{V^{\Gw}(x)}:\abs x=\abs{x_{n}}\right\}=\infty.
\end {equation}
Thus $u$ satisfies (\ref{sing-k10}); \rth {gendom2} and (\ref{sing-k11}) imply that (\ref{sing2}) holds.\qeda\\

The assumption of positivity on $u$ can be weakened if a better a priori estimate is already known. The next result extends \cite [Th 1.2]{FV} into the framework of boundary singularities.
%%%%%%%%%%%%%%%%%%%%%%%%%%%%%%%%%%%%%%%%%%%%%%%%%%%%%
%%%%%%%%%%%%%THEOREM-CLASS 2%%%%%%%%%%%%%%%%%%%%%%%%%%%%%%%%%%%%%%%%%%%%%%%%%%%%%%%%%%%%%%%%%%%%%%%%%%%%%%%%%%%%
\bth {class2}Assume $N-1<q<2N-1$, $\Gw$ is a bounded domain with a $C^2$ boundary, $a\in\prt\Gw$ and $-{\bf e}_{N}$ is the outward normal unit vector to $\prt\Gw$ at $a$. Let $u\in C^1(\overline\Gw\setminus\{a\})$ be a solution of (\ref{equ+abs}) in $\Gw$ vanishing on $\prt\Gw\setminus\{a\}$ such that $u/V^{\Gw}_{a}$ is bounded in $\Gw$. Then there exists $k\in\BBR$ such that  $u=u_{k,a}$.
\es
%%%%%%%%%%%%%%%%%%%%%%%%%%%%%%%%%%%%%%%%%%%%%%%%%%%%%
%%%%%%%%%%%%%PROOF OF THEOREM-CLASS 2%%%%%%%%%%%%%%%%%%%%%%%%%%%%%%%%%%%%%%%%%%%%%%%%%%%%%%%%%%%%%%%%%%%%%%%%%%%%%%%%%%%%
\Proof The outline of the proof are very similar to the finite case of the previous theorem. We still assume $a=0$ and define $k$ by (\ref {S1}). If $k=0$ the maximum principle implies $u\leq 0$ and we return to \rth {class} in the case $u\leq 0$. If $k\neq 0$, $k>0$ for example, (\ref{S2}) and (\ref{S3}) apply. By the previous scaling method we derive that
$u_{n_{k}}$ converges to some function $v$ in the $C^1_{loc}$-topology of $H=\BBR^N_{+}$ which is $N$-harmonic in $H$ and vanishes on $\prt H\setminus\{0\}$. Because $r_{n_{k}}u_{k}(r_{n_{k}}x)$ converges to $kV^{H}_{0}$, the tangency condition of $v$ and 
$kV^{H}_{0}$ at some $\xi$ implies that $v=kV^{H}_{0}$. Thus $u(x)\geq 0$ for $\abs x=r_{n_{k}}$ for $n_{k}$ large enough. This implies that $u\geq 0$ in $\Gw$ and we are back to \rth{class}.\qeda\\

\nind \Remark In the semilinear case of problem (\ref{lin2}), it is proved in \cite {GV} that any signed solution $u$ which satisfies $\lim_{x\to a}\abs {x-a}^{N}u(x)=0$ has constant sign. The exponent $N$ characterize the minimal changing sign harmonic function vanishing on $\prt\Gw$, with an isolated singularity at $a$. Changing sign singular $N$-harmonic functions are constructed in \cite {BV}. In particular there exist singular $N$-harmonic functions $w$ under the form
$$w(r,\gs)=r^{-\gb_2}\gw(\gs)
$$
where 
$$\gb_2=\myfrac{7N-1+\sqrt{N^2+12N+12}}{6(N-1)} $$
and $\gw$ is defined on $S^{N-1}_+=\{x\in S^{N-1}:x_{N}>0\}$, vanishes on the equator $\prt S^{N-1}_+$, is positive on $S^{N-1}_+\cap\{x:x_{N-1}>0\}$ and negative on $S^{N-1}_+\cap\{x:x_{N-1}<0\}$. 
A natural question is therefore wether any signed solution $u$ of (\ref{equ+abs}) in $\Gw$ which vanishes on $\prt\Gw\setminus\{a\}$ and satisfies $\lim_{x\to a}\abs {x-a}^{\gb_2}u(x)=0$ has constant sign, and can be henceforth classified through \rth {class}.\\

\nind {\it Final remark.} If one replaces the $N$-harmonic operator by the $p$-harmonic operator ($p>1$) and tries to extend the results of sections 2, 3, 4, several difficulties will appear. Even if the existence of separable singular solutions is known, the precise value of the exponent 
$\gb>0$ such that $(r,\gs)\mapsto r^{-\gb}\phi(\gs)$ is $p$-harmonic and positive in $H$ and vanishes on $\prt H\setminus \{0\}$ is unknown but for the specific cases $N=2$ or $p=N$ or $p=2$. Notice that in that case the function $\gf$ satisfies the so-called {\it spherical $p$-harmonic spectral equation}
\begin {equation}\label {sfher}\left\{\BA {l}
\!\!\!-div_{\gs}\left(\left(\gb^2\phi^2+\abs {\nabla_{\gs}\phi}^{2}\right)^{(p-2)/2}\!\!\!
\nabla_{\gs}\phi\right)-\gl\left(\gb^2\phi^2+\abs {\nabla_{\gs}\phi}^{2}
\right)^{(p-2)/2}\!\!\!\!\!\phi=0\quad\mbox {in }S^{N-1}_{+}\\
\phantom{;--div_{\gs}\left(\left(\gb^2\phi^2+\abs {\nabla_{\gs}\phi}^{2}\right)^{(p-2)/2}\!\!\!
\nabla_{\gs}\phi\right)-\gl\left(\gb^2\phi^2+\abs {\nabla_{\gs}\phi}^{2}
\right)}
\phi=0\quad\,\mbox { on }\prt S^{N-1}_{+}.\\
\EA\right.\end {equation}
where $\gl=\gb(\gb(p-1)+p-N)$. If $p=2$ then $\gb=N-1$, while if $N=2$, $\gb$ is the positive root of the equation
\begin {equation}\label {KV}
3\gb^2+2\myfrac{p-3}{p-1}\gb-1=0.
\end {equation} 
Furthermore, up to now and due to the lack of conformal invariance, it has not been possible to construct the equivalent of the $V_a^\Gw$ in a general smooth bounded domain $\Gw$, that are positive $p$-harmonic functions in $\Gw$, vanishing on $\prt\Gw\setminus\{a\}$ and satisfying
\begin{equation}\label {weak-p}
 \lim_{\scz {\BA {c}x\to a\\
 (x-a)/\abs {x-a}\to\gs\EA}}\!\!\!\!\!\!\abs {x-a}^{\gb}u(x)=\phi(\gs).
\end{equation}
However, if $\Gw=H=\BBR^N_+$ the removability and the classification results of Sections 3 and 4 are still valid. The proofs of these theorems are developed in \cite{Bo}.
%%%%%%%%%%%%%%%%%%%%%%%%%%%%%%%%%%%%%%%%%%%%%%%
%%%%%%%%%BIBLIOGRAPHY%%%%%%%%%%%%%%%%%%%%%%%%%%%%%
%%%%%%%%%%%%%%%%%%%%%%%%%%%%%%%%%%%%%%%%%%%%%%%%%%%%
\small{
Laboratoire de Math\'ematiques et Physique Th\'eorique\\
CNRS UMR 6083\\
Facult\'e des Sciences\\
Universit\'e Fran\c{c}ois-Rabelais\\
F37200 Tours  France\\

\noindent borghol@univ-tours.fr\\
veronl@lmpt.univ-tours.fr}
%Furthermore%Furthermore%Furthermore
%Furthermore%Furthermore%Furthermore
%Furthermore%Furthermore%Furthermore
 %%END DOCUMENT%%%%%%%%%%%%%%%%%%%%%%%%
\end {document}